\documentclass[a4paper,12pt]{article}

\usepackage[english]{babel}
\usepackage[T1]{fontenc}
\usepackage[utf8x]{inputenc}

\usepackage{verbatim}
\usepackage{amsfonts,amsmath,amssymb,bbding,bbm,color,enumerate,textcomp,url,yfonts,geometry}
\usepackage{blindtext}
\usepackage{multicol}
\usepackage{lipsum}
\usepackage{mathtools}
\usepackage{cuted}

\usepackage{authblk}
\usepackage[export]{adjustbox}

\usepackage{amsmath}
\usepackage{bm}
\usepackage{array}
\usepackage{amssymb}
\usepackage{amscd}
\usepackage{multirow}
\usepackage{calrsfs}
\usepackage{url}
\usepackage{color}
\usepackage{soul}
\usepackage{xcolor}
\usepackage[]{algorithm2e}
\usepackage[normalem]{ulem}
\usepackage{mathrsfs}  
\usepackage{mathtools}
\usepackage{dirtytalk}
\usepackage{enumerate}

\def\ml{l\kern-0.035cm\char39\kern-0.03cm}
\def\mL{L\kern-0.08cm\char39}

\definecolor{brown}{rgb}{0.43, 0.21, 0.1}
\newcommand{\js}[1]{{\color{brown}#1}}
\newcommand{\jb}[1]{{\color{red}#1}}

\newtheorem {definicia}{Definition}[section]
\newtheorem {veta}[definicia]{Theorem}
\newtheorem {lema}[definicia]{Lemma}
\newtheorem {tvrdenie}[definicia]{Proposition}

\newtheorem {example}[definicia]{Example}
\newtheorem {remark}[definicia]{Remark}
\newtheorem {dosledok}[definicia]{Corollary}

\newcommand{\dokaz}{\medskip \noindent {\bf Proof. \ \ }}
\newcommand{\dokazp}[1]{\medskip \noindent {\bf Proof of #1.\ \ }}

\newcommand{\Borel}{\mathcal{B}}

\newcommand{\counting}{\nu_{\#}}

\newcommand{\x}{\mathbf{x}}

\newcommand{\bin}{\mathbf{1}}

\newcommand{\N}{\mathbb{N}}
\newcommand{\R}{\mathbb{R}}
\newcommand\sA{\mathsf{A}}
\newcommand\aA[2][E]{\sA(#2|#1)}
\newcommand\aAi[3][E]{\sA^{#3}(#2|#1)}

\newcommand\cA{{\mathcal A}} 

\newcommand\bA{{\mathscr {A}}}

\newcommand\cE{{\mathcal{E}}}

\newcommand\cS{{\mathcal{S}}}

\newcommand\hE{\hat{E}}

\newcommand\Chint{\mathrm{C}_{\cA}(f,\mu)}

\newcommand{\proof}{\medskip \noindent {\bf Proof. \ \ }}
\newcommand{\collectionGSF}{\mathcal{E}}

\def\gsf{\mu_{\bA}(\x,\alpha)}

\def\gsf#1{\mu_{\mathcal{A}^{#1}}(\x,\alpha)}

\def\gsf#1#2#3{\mu_{\mathcal{A}^{#1}}(#2,#3)}

\newcommand{\qed}{\null\hfill $\Box\;\;$ \medskip}

\usepackage{cite}

\begin{document}

\begin{center}
{\bf {\sc \Large Conditional Aggregation-Based Choquet Integral as~a~Choquet Integral on a~Hyperspace}}
\end{center}

\vskip 12pt

\begin{center}
{\bf Jana Borzov\'{a}, Lenka Hal\v{c}inov\'{a}, Jaroslav \v{S}upina}\footnote{
\newline Corresponding author: jana.borzova@upjs.sk}
\end{center}

\vskip 20pt

\begin{abstract}
We aim at representing the~recently introduced conditional aggregation-based Choquet integral as a~standard Choquet integral on a~hyperset. The~representation is one of transformations considered by R.R.~Yager and R.~Mesiar in~2015. Thus we study the~properties of the~conditional aggregation-based Choquet integral using well-known facts about the~standard Choquet integral. In~particular, we obtain several formulas for its computation. We also provide the representation of the~conditional aggregation-based Choquet integral in terms of the~M{\"o}bius transform. 
\end{abstract}

\section{Introduction}{

We study the concept of generalized survival functions and the corresponding conditional aggregation-based Choquet integral based on them. Both notions were recently introduced in~\cite{BoczekHalcinovaHutnikKaluszka2020}, and were inspired by a~concept of size introduced in~\cite{DoThiele2015}, and further developed in~\cite{HalcinovaHutnikKiselakSupina2019}. The~main ingredient of this construction are conditional aggregation operators that cover many existing aggregations, such as the~arithmetic and geometric mean, or plenty of~integrals known in~the~literature~\cite{KlementMesiarPap2010, Shilkret1971, Sugeno1974}. Despite some similarity of the definition of conditional aggregation operator with that of the interaction operator~\cite{HondaOkazaki2017}, the authors in~\cite{BoczekHutnikKaluszka2022} showed that the concept of conditional aggregation operators is more general than that of extended interaction operators. The~new concept of integrals offers many variations depending on chosen aggregation operators.
The generalized survival function and the corresponding Choquet integral are interesting from the point of view of some applications.
Using the new concept, the authors model the customer satisfaction with the purchase or problems in insurance, see~\cite[Example 6.1--Example 6.3]{BoczekHalcinovaHutnikKaluszka2020}. In~\cite{BasarikBorzovaHalcinovaSupina2023} the authors used this concept to solve the Knapsack problem and the problem of accommodation options.


The main contribution of the paper is to provide a different perspective on the construction of the~conditional aggregation-based Choquet integral compared to the original paper~\cite{BoczekHalcinovaHutnikKaluszka2020}. We show that the new Choquet integral is a standard Choquet integral on a hyperset.

As in the classic case of the Choquet integral, we consider nonnegative real-valued functions that are defined on a nonempty set $X$. Further, we consider a monotone measure $\mu$ defined on a collection of subsets of $X$~denoted~by~$\hat{\cE}$, i.e. $\hat{\cE}\subseteq 2^X$.

The construction of the~conditional aggregation-based Choquet integral does not use level sets of real-valued functions. A~conditional aggregation operator $\aA[E]{\cdot}$ assigns to $f\colon X\to[0,+\infty]$ a new value $\aA[E]{f}$ from $[0,+\infty]$. The~assigned value depends on a~conditional set~$E$. In fact, conditional sets are complements of sets from $\hat{\cE}$. 

The~family~$\cA$ of operators $\aA[E]{\cdot}$ for conditional sets~$E$ (which are determined by the collection $\hat{\cE}$ selection), together with a~monotone measure~$\mu$, is then used to build a~generalized survival function $\gsf{}{f}{\alpha}$.
The~generalized survival function becomes a~key to introduce conditional aggregation-based Choquet integral of a~function~$f$ as 
$$\mathrm{C}_{\cA}(f,\mu)=\int_{0}^{\infty}\gsf{}{f}{\alpha}\,\mathrm{d}\alpha.$$
We show that a~monotone measure space~$(X,\mu)$ may be assigned a~monotone measure space~
$(\hat{\cE},\mathtt{N}_\mu)$ with~$\mathtt{N}_\mu$ introduced in Section~\ref{measure_to_measure}, and a~function~$f$ on~$X$ may be transformed to a~function~$T_{\cA,f}$ on~
$\hat{\cE}$, such that $$\mathrm{C}_{\cA}(f,\mu)=\mathrm{Ch}(T_{\cA,f},\mathtt{N}_{\mu}).$$ 
By $\mathrm{Ch}(T_{\cA,f},\mathtt{N}_{\mu})$ we mean the~standard Choquet integral of~a~function~$T_{\cA,f}$ computed on a~measure space~
$(\hat{\cE},\mathtt{N}_\mu)$, i.e., $$\mathrm{Ch}(T_{\cA,f},\mathtt{N}_{\mu})=\int_{0}^{\infty}\mathtt{N}_{\mu}(\{T_{\cA,f}>\alpha\})\,\mathrm{d}\alpha.$$

Actually,  the monotone measure $\mathtt{N}_\mu$ is one of transformations of~monotone measures to~the~power set, studied by~Yager and 
Mesiar in~\cite{YagerMesiar2015}. Indeed, the domain of monotone measure~$\mathtt{N}_\mu$ is $2^{\hat{\cE}}$.

Such a grasp of~the~concept can be beneficial from several points of view:
\begin{enumerate}[\rm (i)]
    \item The~generalized Choquet integral may be seen as~a~standard Choquet integral based on $\mathtt{N}_\mu$-survival function (Corollary~\ref{general=classic}). In particular, this is the~case of size-based integrals as of~\cite{DoThiele2015,HalcinovaHutnikKiselakSupina2019}, see Remark~\ref{size-integral}. 
    \item The properties of fuzzy integrals based on generalized survival function, of which the Choquet integral is most often mentioned, see~\cite{BoczekHalcinovaHutnikKaluszka2020,Halcinova2017}, strongly depend on~properties of~$\mathtt{N}_\mu$. Thus the study of its properties is reasonable, for more details see Section~\ref{Choquet_integral}. For some properties of~$\mathtt{N}_\mu$ we refer to Section~\ref{measure_to_measure}, for particular result where the~property of $\mathtt{N}_\mu$ was used to derive the~property of~the~generalized Choquet integral we refer to~Proposition~\ref{linearity}. We show a~connection to~integrals studied in~\cite{DimuroFernandezBedregal2020, LuccaDimuro2019}, see~Remark~\ref{fusion}.
    \item The~application of the~standard M{\"o}bius transform to~monotone measure $\mathtt{N}_\mu$ leads to~a~new formula for~computation the~generalized Choquet integral, see Section~\ref{Mobius_sec}.
    We derive also other formulas for the~generalized Choquet integral, see Proposition~\ref{discreteChoquet}.
\end{enumerate}}

\medskip

The~paper is organized as follows. We recall some terminology used throughout the~paper at the~end of this introduction. Since the~terminology on minitive monotone measure is not unified in the~literature, we devote a~separate Section~\ref{S-minitive} to set up one we shall use, and to recall some standard results on it. Section~\ref{measure_to_measure} introduces the~main tool of the~paper, the~particular transformation of one monotone measure space to the~other. We explore basic properties of the~transformed monotone measure, and we illustrate the~transformation on many examples, so the~reader could get familiarised with the~monotone measure on a~hyperset. In Section~\ref{S-aggr+surv} we present our main result in Theorem~\ref{Vtransform}. It states that the~value of the~generalized survival function~$\gsf{}{f}{\alpha}$ of function~$f$ on~$X$ in~$\alpha$ is equal to $\mathtt{N}_\mu$-measure of $\alpha$-cut of a~function~$T_{\cA,f}$ on~$\hat{\cE}$. A~brief summary on conditional aggregation operators, and the~generalized survival function is included in~the~same section. 

In Section~\ref{Choquet_integral} we provide a new perspective on the conditional aggregation-based Choquet integral construction. It is included in Corollary~\ref{general=classic} which follows directly from Theorem~\ref{Vtransform}. Representation of the new Choquet integral in terms of the M{\"o}bius transform is derived in Section~\ref{Mobius_sec}. In the~last section, we compute a~dual monotone measure to monotone measure~$\mathtt{N}_\mu$, and once again apply monotone measure~$\mathtt{N}_\mu$ in an~alternative proof of a~result by~M.~Boczek et al.~\cite{BoczekHutnikKaluzskaKleinova2021}. Finally, we include the~Conclusion, and Appendix with a~couple of~proofs.


In the paper we shall use the standard terminology. Since we work with two measure spaces $(X,\mu)$ and $(\hat{\cE},\mathtt{N}_\mu)$, let us introduce a common notion $\Omega$ for a nonempty basic set to~cover both cases $X$ and $\hat{\cE}$, respectively. Let $\cS$ be a nonempty collection of subsets of $\Omega$.
A~set function $m\colon\mathcal{S}\to[0,\infty]$  will be called a \textit{monotone set function} iff $m$ is~nondecreasing, i.e., $m(E)\leq m(F)$ whenever $E\subseteq F$, $E,F\in\mathcal{S}$. If $\emptyset\in\mathcal{S}$, then, in fact, the pair $(\Omega,\mathcal{S})$ is a prespace according to {\v{S}}ipo{\v{s}}, see~\cite{Sipos1979,MesiarSipos1994}. Naturally, we consider $m(\emptyset)=0$ and then we use a notion \textit{a monotone measure} instead of a monotone set function. In this sense, a monotone measure is identical to the {\v{S}}ipo{\v{s}} premeasure.
If $\Omega\in\cS$ and $m(\Omega)=1$, then a~monotone measure~$m$ will be called \textit{a capacity}. 
Moreover, we shall suppose $m(\Omega)>0$. We shall use short notation to denote cuts of various nonnegative real-valued functions, including monotone set functions. Indeed, if $H\colon\Omega\to[0,+\infty]$,
\[
\{H>\alpha\}=\{x\in\Omega: H(x)>\alpha\},
\]
for each $\alpha\in\R$, and similarly for other symbols among $=,\geq,<,\leq$. The discrete set $\{1,\dots,n\}$, $n\geq 1$, is denoted by $[n]$.

\section{Minitive monotone measures}\label{S-minitive}

Before introducing the~fundamental tool for this section let us recall the~terminology on~maxitive and minitive measures that are widely studied in the literature, see e.g.~\cite{KRATSCHMER2003455,Poncet2017,Puhalskii,Shilkret1971,WangKlir2009}.
The maxitive measures were coined already 
by Shilkret in~\cite{Shilkret1971} and have been widely used, especially in the fields of probability theory and fuzzy
theory,~\cite{Poncet2017}. One can find many other terms in the literature for maxitive measures, e.g. $\tau$-maxitive measures~\cite{Puhalskii}, supremum-preserving measures~\cite{KRATSCHMER2003455}, completely maxitive measures~\cite{Pap1995NullAdditiveSF,Poncet2017}. Usually, 
if the capacity is maxitive, then it is called a \textit{possibility measure}~\cite{Dubois2015,WangKlir2009,ZADEH19783}. The conjugate of a maxitive monotone measure is the minitive monotone measure. If the capacity is minitive, then it is called a \textit{necessity measure}.


\begin{definicia}
Let $m\colon \mathcal{S}\to[0,\infty]$ with $\emptyset\in\mathcal{S}\subseteq2^{\Omega}$ be a~monotone measure.
Monotone measure $m$ is called
\begin{enumerate}[\rm (a)]
    \item  \textit{minitive} iff
    \[
m\left(\bigcap\limits_{t\in T}E_t\right)=\inf\limits_{t\in T}m(E_t)
\]
for any subfamily $\{E_t\colon t\in T\}$ \text{of}~$\mathcal{S}$ whose intersection is in $\mathcal{S}$ with $T$ being an arbitrary~index set.
\item  \textit{maxitive} iff
    \[
m\left(\bigcup\limits_{t\in T}E_t\right)=\sup\limits_{t\in T}m(E_t)
\]
for any subfamily $\{E_t\colon t\in T\}$ \text{of}~$\mathcal{S}$ whose union is in $\mathcal{S}$ with $T$ being an arbitrary~index set.
\end{enumerate}
\end{definicia}

We recall a~generic way to construct minitive and maxitive measures that are known in~the~literature, e.g. see this result in~\cite[Theorem 4.23]{WangKlir2009} for possibility measures. The~same is true for maxitive, minitive measures, as well.
The~generating function is in the literature called distribution of $m$ (possibility distribution for possibility measures in~\cite{Dubois2015}) or possibility profile in~\cite{WangKlir2009}. We include a~proof of Proposition~\ref{necesschar} in the~Appendix.

\begin{tvrdenie}\label{necesschar}
Let $\Omega$ be a~nonempty set, and let $m\colon 2^\Omega\to[0,\infty]$. 
\begin{enumerate}[\rm (a)]
    \item $m$ is minitive monotone measure if and only if there is $\pi'\colon\Omega\to[0,\infty]$ such that 
\[
m(E)=\inf\{\pi'(x): x\in\Omega\setminus E\}\footnote{Infimum is computed with respect to the~set $\{\pi'(x): x\in\Omega\}$, i.e., $\inf\emptyset=\max\{\pi'(x): x\in\Omega\}$ if it exists, otherwise $\inf\emptyset=\infty$. }
\]
for any $E\subseteq\Omega$.
    \item $m$ is maxitive measure if and only if there is $\pi\colon\Omega\to[0,\infty]$ such that 
\[
m(E)=\sup\{\pi(x): x\in E\}
\]
for any $E\subseteq\Omega$.
\end{enumerate}
\end{tvrdenie}
\par

In the following we put stress on minitive measures, since they will play a crucial role in~this paper. We shall denote them by $N$ (on the other hand, the maxitive measures we shall denote by $\Pi$). Let us summarize their basic properties, see similar results for possibility measures in~\cite{WangKlir2009}. A monotone set function $m:\mathcal{S}\to[0,\infty]$ will be called
\begin{enumerate}[\rm (i)]
\item \textit{continuous from above} or \textit{upper semicontinuous} iff $\lim\limits_{n\to\infty} m(E_n)=m\left(\mathop{\cap}\limits_{n=1}^{\infty}E_n\right)$ for any $(E_n){_1^\infty}\subseteq\mathcal{S}$, $E_1\supseteq E_2
\supseteq\dots$, $\mathop{\cap}\limits_{n=1}^{\infty}E_n\in\mathcal{S}$;
    \item \textit{continuous from below} or \textit{lower semicontinuous} iff $\lim\limits_{n\to\infty} m(E_n)=m\left(\mathop{\cup}\limits_{n=1}^{\infty}E_n\right)$ for any $(E_n){_1^\infty}\subseteq\mathcal{S}$, $E_1\subseteq E_2
\subseteq\dots$, $\mathop{\cup}\limits_{n=1}^{\infty}E_n\in\mathcal{S}$;
\item \textit{continuous} if it is continuous from below and above;
\item \textit{superadditive} iff $m(E\cup F)\geq m(E)+m(F)$ whenever $E\cup F\in\mathcal{S}$, $E, F\in\mathcal{S}$, $E\cap F=\emptyset$.
\end{enumerate}
\smallskip

\begin{tvrdenie}\label{necessvlastnosti}
Let $N:2^\Omega\to[0,\infty]$ 
be~a~minitive monotone measure induced by a~function $\pi'\colon \Omega\to[0,\infty]$.
Then 
\begin{enumerate}[\rm (a)]
    \item $N$ is superadditive and continuous from above.
      \item $N(\Omega\setminus\{x\})=\pi'(x)$.
    \item If $\pi'(x)=0$ and $x\not\in E$ then $N(E)=0$. 
    \item $N(\Omega)=\sup\{\pi'(x): x\in \Omega\}$. 
    \item $N(E)=\inf\{N(\Omega\setminus\{x\}): E\subseteq \Omega\setminus\{x\}\}$. 
\end{enumerate}
\end{tvrdenie}
\par

\begin{remark}
Once a~minitive {measure}~$N$ has just finite values, we can define its dual maxitive {measure}~$\Pi$ by $\Pi(E)=N(\Omega)-N(\Omega\setminus E)$ and vice versa. 
 Then associated function $\pi'$ is bounded and we can define a~dual object, a~function $\pi$ by $\pi(x)=\max\{\pi'(y): y\in \Omega\}-\pi'(x)$. Function $\pi$ induces a~dual maxitive measure~$\Pi$ to a given $N$. \end{remark}

\section{Measure-to-measure transform }\label{measure_to_measure}

As we have stated in Introduction, since we are mainly interested in (the properties of) generalized survival functions and integrals based on them, we focus on the properties of~minitive measures, in particular on monotone measure $\mathtt{N}_\mu$ below.  
Let us consider a~nonempty set~$X$, $\hat{\cE}\subseteq2^X$, and a~monotone set function~$\mu:\hat{\cE}\to[0,\infty]$. The class of all monotone set functions on $(X,\hat{\cE})$ we shall denote by  $\mathbf{M}$. The elements of $2^{\hat{\cE}}$ will be denoted by $\hat{E}$, i.e. $\hat{E}\in2^{\hat{\cE}}$ or $\hat{E}\subseteq\hat{\mathcal{E}}$.  

Now, we shall introduce a~minitive monotone measure\footnote{The domain of $\mathtt{N}_\mu$ is a~subset of $2^{2^X}$.}  
$\mathtt{N}_\mu: 2^{\hat{\cE}}\to[0,\infty]$  tied with~a~monotone set function~$\mu\in\mathbf{M}$:
\begin{align}\label{def:Nmu}
\mathtt{N}_\mu(\hE):=\inf\{\mu(E)\colon E\in\hat{\cE}\setminus\hE\}.\end{align}
with the convention  $\inf\emptyset=\mu(X)$ whenever $X\in\hat{\cE}$.\footnote{This follows from convention in Section 2.} Let us point out that $\mathtt{N}_\mu$ is defined on~a~powerset of~$\hat{\cE}$ not only on~a~subset of~$2^X$. It is easy to see that $\mathtt{N}_\mu$ is a~minitive monotone measure. The~role of function $\pi'$ described in~Proposition~\ref{necesschar} is played by a monotone set function $\mu$, and $\Omega=\hat{\cE}$.
For a~better comprehension of the introduction of~the~monotone measure $\mathtt{N}_{\mu}$, let us give several easy examples.

\par
\begin{example}[Counting measure]\label{countingmeasure} Let us consider 
$\hat{\cE}=2^{[3]}$ and a counting measure $\counting\colon 2^{[3]}\to \{0,1,2,3\}.$ Now, $\mathtt{N}_{\counting}(\hE)=\min\{|E|\colon E\in2^{[3]}\setminus\hE\}$ for each $\hat{E}\subseteq 2^{[3]}$. Then we have
{{\begin{align*}\mathtt{N}_{\counting}(\hE)=\begin{cases}0\quad&\textrm{if}\,\,\, \emptyset\notin\hE
;\\ 1\quad&\textrm{if}\,\,\,\left(\exists E\notin\hE\right)|E|=1\text{ and}\,\, \emptyset\in\hE;\\ 
2\quad&\textrm{if}\,\,\,\left(\exists E\notin\hE\right)|E|=2\text{ and}\,\,\hE\supseteq\{F\colon|F|< 2\};\\ 
3\quad&\textrm{if}\,\,\,\hE=2^{[3]}, 2^{[3]}\setminus \{[3]\}
.\ \end{cases}\end{align*}}}
\end{example}
\par

If a~monotone set function~$\mu$ achieves the~finite number of~values, then the corresponding measure $\mathtt{N}_{\mu}$ may be expressed in the~simpler form. In what follows, let us recall that 
\begin{align*}
\{\mu<\alpha\}&=\{E\in\hat{\cE}:\,\,\mu(E)<\alpha\},\\
\{\mu\leq\alpha\}&=\{E\in\hat{\cE}:\,\,\mu(E)\leq\alpha\}.
\end{align*}
Then from the definition we have 
\[
\mathtt{N}_\mu(\{\mu\leq\alpha\})\geq\mathtt{N}_\mu(\{\mu<\alpha\})\geq\alpha.
\]

\begin{example}[Finitely many values]\label{konecnypocet} Let us consider 
$(u_j)_1^q$ be an increasing enumeration of all values of a~monotone set function~$\mu:\hat{\cE}\to[0,\infty]$, i.e.,
$$u_1<u_2<\dots<u_q.$$ Then the~transformation $\mathtt{N}_\mu: 2^{\hat{\cE}}\to\{u_1,\dots,u_q\}{\cup \{0\}}$ is a~monotone measure given for each~$\hE\subseteq \hat{\cE}$ as 
{\begin{align*}
\mathtt{N}_{\mu}(\hE)=\begin{cases}0\quad&\textrm{if}\,\,\,\hat{E}=\emptyset;\\
u_j&\textrm{if}\,\,\,\left(\exists E\notin\hE\right)\mu(E)= u_j
\,\, \text{and}\,\, \hE\supseteq\{\mu<u_j\};\\
u_q\quad&\textrm{if}\,\,\,\hat{E}={\hat{\mathcal{E}}}
.\ \end{cases}\end{align*}}
\end{example}

Notice that, although $\emptyset\notin\hat{\cE}$, $\mathtt{N}_{\mu}$ always achieves the zero value because of its definition. Clearly, a monotone set function $\mu$ from Example~\ref{konecnypocet} covers the counting measure. The~fact that $\mathtt{N}_{\counting}(\emptyset)=0$ follows from~the first case of~the~formula in Example~\ref{countingmeasure}. Taking $j\in\{1,2,3\}$ and $q=3$ in Example~\ref{konecnypocet}, we derive the other cases in formula in Example~\ref{countingmeasure}.

As a~special case of preceding example is also $\{0,1\}$-valued monotone measure. In particular, the~transformation $\mathtt{N}_\mu$ is also a $\{0,1\}$-valued measure. We can compute the~transformation~$\mathtt{N}_\mu$ for the~weakest and the~strongest measures. 

\begin{example}[Weak \& strong]
Let $\{\emptyset,X\}\subseteq\hat{\cE}$. For~the~weakest capacity $\mu_{*}:{\hat{\cE}}\to\{0,1\}$ (i.e., $\mu_{*}(E)=1$ for $E=X$ and $\mu_{*}(E)=0$ otherwise) the~transformation $\mathtt{N}_{\mu_{*}}:2^{\hat{\cE}}\to\{0,1\}$ takes the form 
{{\[
\mathtt{N}_{\mu_{*}}(\hat{E})=
\begin{cases}
1&\textrm{if}\,\, \hat{E}=\hat{\cE},\, \hat{\cE}\setminus \{X\};\\
0&\text{otherwise}.
\end{cases}
\]}}
On the other hand, for the strongest capacity $\mu^{*}:\hat{\cE}\to\{0,1\}$ (i.e., $\mu^{*}(E)=0$ for $E= \emptyset$ and $\mu^{*}(E)=1$ otherwise) the~transformation $\mathtt{N}_{\mu^{*}}:2^{\hat{\cE}}\to\{0,1\}$ takes the~form 
{{\[
\mathtt{N}_{\mu^{*}}(\hat{E})=
\begin{cases}
0&\textrm{if}\,\,\hE=\emptyset,\, \hat{\cE}\setminus\{\emptyset\};\\
1&\text{otherwise}.
\end{cases}
\]}}

\end{example}

\begin{example}[Possibility \& necessity] Let us consider a~finite set $X=[n]$ and $\hat{\cE}=2^{[n]}$. We fix a~permutation~$(\cdot)$ of elements of $[n]$ such that $0=\pi((0))\leq \pi((1))\leq \pi((2))\leq\dots\leq \pi((n))=1$, and we set $E_{(i)}=\{(i),\dots,(n)\}$ for any $i\in[n]$ with the convention $E_{(n+1)}=\emptyset$. By $\tau$ we denote a~system of elements of~$[n]$ defined as 
\begin{equation}\tau:=\{i\in[n]\colon \pi((i-1))<\pi((i))\}.\footnote{If we have $X=[4]$ and $\pi(1)=0.7$, $\pi(2)=0.4$, $\pi(3)=1$ and $\pi(4)=0.4$, then we can take the permutation $(1)=2$, $(2)=4$, $(3)=1$, $(4)=3$. Finally, $\tau=\{1,3,4\}$.}
\end{equation}
\begin{enumerate}[\rm(a)]
    \item Let $\Pi\colon 2^{[n]}\to[0,1]$ be a~possibility measure, i.e., $\Pi(E)=\max\{\pi(i)\colon i\in E\}$.
  Then the monotone measure $\mathtt{N}_\Pi$ takes the form 
{
\begin{align*}\mathtt{N}_{\Pi}(\hE)=\begin{cases}0&{\textrm{if}\,\,\emptyset\notin\hE};\\
\pi((i))&\textrm{if}\,\,i\in\tau, \left(\exists E\notin\hE\right)\Pi(E)=\pi((i))\,\, \text{and}\,\,  \hE\supseteq2^{E_{(i)}^c};
\\
1\quad&\textrm{if}\,\,\,\hE=2^{[n]}.
\end{cases}\end{align*}
}
Observe that for every $i\in\tau$ it holds that $\pi((i))>0$. The condition for taking the value $\pi((i))$, $i\in\tau$, we directly obtain from Example~\ref{konecnypocet}. Indeed, it has to hold $\hE\supseteq\left\{F:\max_{(k)\in F}\pi((k))<\pi((i))\right\}$. 

\item Let $N\colon 2^{[n]}\to[0,1]$ be a~necessity measure, i.e., $N(E)=1-\Pi(E^c)$. Let us define $k=\max\{j\in[n]: \pi((j))=\pi((i))\}$. Then the monotone measure $\mathtt{N}_N$ takes the form {
\begin{align*}\mathtt{N}_{N}(\hE)=\begin{cases}0\quad&\textrm{if}\,\,\hat{E}=\emptyset;\\
1-\pi((i))\quad&\textrm{if}\,\,i\in\tau,
\left(\exists E\notin\hE\right)\Pi(E^c)=\pi((i)) \,\, \text{and}\,\,\left(\forall E\notin\hE\right) E\supseteq E_{(k+1)};
\\1\quad&\textrm{if}\,\,\,\left(\forall E\notin\hE\right)\Pi(E^c)=0.
\end{cases}\end{align*} }

\end{enumerate}
\end{example}

All previous examples were based on $\hat{\cE}$ being finite. For~infinite~$\hat{\cE}$  we refer to Example~\ref{Lebesgue}. Let us summarize the~very basic properties of $\mathtt{N}_\mu\colon 2^{\hat{\cE}}\to[0,\infty]$, mostly consequences of~Proposition~\ref{necessvlastnosti}. Hence, $\mathtt{N}_\mu$ is a~monotone measure which is in~addition minitive, superadditive and upper semicontinuous. The~family of all sets in~{$2^X$} with $\mu$-measure zero is denoted by~$\mathcal{N}_{\mu}$.
\begin{tvrdenie}\label{vlastnosti}
Let $\mu\in\mathbf{M}$. 
\begin{enumerate}[\rm (a)]
    \item If $\mu(E)=0$ and $E\not\in\hE$, then $\mathtt{N}_\mu(\hE)=0$. 
    \item If $\emptyset\in\hat{\cE}$ and $\emptyset\not\in\hE$, then $\mathtt{N}_\mu(\hE)=0$.
    \item If $\emptyset\in\hat{\cE}$, and $\mu(E)\neq0$, then $\mathtt{N}_{\mu}(\{E\})=0$.
    \item If $\emptyset\in\hat{\cE}$ and $|\mathcal{N}_{\mu}|>1$, then $\mathtt{N}_{\mu}(\{\emptyset\})=0$.
    \item $\mathtt{N}_\mu(\hE)=0$ if  for every $\varepsilon>0$ there is $E\in\hat{\cE}\setminus\hE$ such that $\mu(E)<\varepsilon$.
    \item If $X\in\hat{\cE}$, then $\mathtt{N}_\mu(\hat{\cE})=\mu(X)$. 
    \item If $E\in\hat{\cE}$, then $\mathtt{N}_\mu(\hat{\cE}\setminus\{E\})= \mu(E)$.
\end{enumerate}
\end{tvrdenie}

\par Assuming that $\hat{\cE}$ is finite leads to some improvements of the~preceding assertion. In fact, we obtain a~complete characterization of $\mathtt{N}_\mu$-measure zero sets.
\begin{lema}
Let $\hat{\cE}$ be finite such that $\emptyset\in\hat{\cE}$ and $\mu\in\mathbf{M}$.
\begin{enumerate}[\rm (a)]
    \item $\mathtt{N}_\mu(\hE)=0$ if and only if there is $E\not\in\hE$ such that~$\mu(E)~=~0$.
    \item $\mathtt{N}_{\mu}(\{\emptyset\})=0$ if and only if $|\mathcal{N}_\mu|>1$.
    \item $\mathtt{rng}(\mathtt{N}_{\mu})=\mathtt{rng}(\mu)$.
\end{enumerate}
\end{lema}
\dokaz
One implication of~(a) follows from Proposition~\ref{vlastnosti}, and the other one by the fact that if the~minimum of a~set is zero then the~set has to contain zero number. (b) is a~direct consequence of (a). Part (c) follows from the~definition of $\mathtt{N}_{\mu}$.
\qed
\medskip
\par

It is obvious that monotone measure $\mathtt{N}_\mu$ need not be additive in general. In~fact, by the following~Proposition~\ref{additivitymmi1} only very specific monotone measure $\mu$ induces additivity of $\mathtt{N}_\mu$. \begin{tvrdenie}\label{additivitymmi1}
Let $\{\emptyset,X\}\subseteq\hat{\cE}$, and $\mu\in\mathbf{M}$ with $\mu(X)=c>0$. Then the~following assertions are equivalent:
\begin{itemize}
\item[\rm (i)]
$\mu(E)=0$, if $E=\emptyset$ and $\mu(E)=c$, otherwise;
\item[\rm (ii)]$\mathtt{N}_\mu$ is additive.
\item[\rm (iii)]$\mathtt{N}_\mu$ is $\sigma$-additive.
\end{itemize}
\end{tvrdenie}

\dokaz
Observe that the implication $(iii)\Rightarrow (ii)$ is trivial.

$(i)\Rightarrow (iii)$
Let us consider a sequence $\{\hE_n\}_{n=1}^\infty$ of pairwise disjoint sets in $2^{\hat{\cE}}$ with $\emptyset\notin\hE_n$ for any $n\geq 1$. Then for each $n\geq 1$ it holds that $\mathtt{N}_\mu(\hE_n)=0$ and also $\mathtt{N}_\mu\left(\bigcup\limits_{n=1}^\infty\hE_n\right)=0$. 

Without loss of generality let us suppose $\emptyset\in\hE_1$. Then $\mathtt{N}_\mu(\hE_1)=c$ and for each $n\geq 2$ we have $\mathtt{N}_\mu(\hE_n)=0.$ Consequently, $\mathtt{N}_\mu\left(\bigcup\limits_{n=1}^\infty\hE_n\right)=c$.

$(ii)\Rightarrow (i)$
Suppose that $\mathtt{N}_\mu$ is additive, i.e., for all $\hE_1, \hE_2\in\hat{\cE}$ with $\hE_1\cap\hE_2=\emptyset$ it holds that $\mathtt{N}_\mu(\hE_1)+ \mathtt{N}_\mu(\hE_2)= \mathtt{N}_\mu(\hE_1\cup\hE_2)$. 
Set $\hE_1=\{\emptyset, X\}$ and $\hE_2=\{\hat{\cE}\setminus\{\emptyset, X\}\}$. Then it is clear that $\mathtt{N}_\mu(\hE_2)=0$ and $\mathtt{N}_\mu(\hE_1\cup\hE_2)=\mu(X)=c$ and the~additivity of~$\mathtt{N}_{\mu}$ implies $\mathtt{N}_\mu(\hE_1)=c$. Moreover, it means that for each $E\in\hat{\cE}\setminus\{\emptyset, X\}$ it holds that $\mu(E)=c$. 
\qed \medskip 

In the~realm of~Proposition~\ref{additivitymmi1} one may ask whether there is any other relation between $\mathtt{N}_\mu$-measures of sets and their union. However, the~following fact has to be taken into account. If $E\notin\mathcal{N}_\mu$, then $\mathtt{N}_\mu(\{E\})=0$. Moreover, by Proposition~\ref{vlastnosti}, (g), we have $\mathtt{N}_\mu(\hat{\cE}\setminus\{E\})=\mu(E)$. Finally, $\mathtt{N}_\mu((\hat{\cE}\setminus\{E\})\cup\{E\})=\mu(X)$ by Proposition~\ref{vlastnosti}, (f). Note that $\mu(E)$ may be an~arbitrary positive value of a~monotone measure~$\mu$.

The~cuts of measure~$\mu$ defined earlier characterize the~behavior of monotone measure~$\mathtt{N}_\mu$.

\begin{tvrdenie}\label{cuts}
Let $\mu\in\mathbf{M}$. 
\begin{enumerate}[\rm (a)]
    \item $\alpha\leq\mathtt{N}_\mu(\hat{E})$ if and only if $\{\mu<\alpha\}\subseteq\hat{E}$. 
    \item $\mathtt{N}_\mu(\hat{E})<\beta$ if and only if there is $E\in\{\mu<\beta\}\setminus\hat{E}$.
    \item $\mathtt{N}_\mu(\hat{E})=\alpha$ if and only if $\{\mu<\alpha\}\subseteq\hat{E}$ and for any $\varepsilon>0$ there is $E\in\{\mu<\alpha+\varepsilon\}\setminus\hat{E}$. 
\end{enumerate}
\end{tvrdenie}

Once measure~$\mu$ attains all the~values of~the~interval, the~cuts of monotone measure~$\mu$ describe completely the~behavior of~monotone measure~$\mathtt{N}_\mu$.
In the following example ${\Borel}$ denotes the~Borel $\sigma$ - algebra.
\begin{example}[Lebesgue measure]\label{Lebesgue}
Let us consider $X=[0,1]$ and the~Lebesgue measure~$\lambda$. One can see that for~the~monotone measure $\mathtt{N}_\lambda\colon 2^
{\Borel}\to [0,1]$ we have
\[
\mathtt{N}_\lambda(\{\lambda\leq\alpha\})=\mathtt{N}_\lambda(\{\lambda<\alpha\})=\alpha.
\]
We list the~values of~$\mathtt{N}_\lambda$ for some sample sets.
\begin{enumerate}[\rm (a)]
    \item $\mathtt{N}_\lambda(\mathcal{N}_\lambda)=0$. 
    \item If $\mathcal{N}_\lambda\setminus\hat{E}\neq\emptyset$ then $\mathtt{N}_\lambda(\hat{E})=0$. In particular, if $\emptyset\notin \hat{E}$ or $\{b\}\notin \hat{E}$ then $\mathtt{N}_{\lambda}(\hat{E})=0$.
    \item If $\hat{E}$ contains the finite number of sets, then $\mathtt{N}_\lambda(\hat{E})=0$.
    \item $\mathtt{N}_\lambda(\Borel\setminus\left\{\left(b-\frac{1}{n},b+\frac{1}{n}\right),\,\, n\in\N\}\right\})=0.$
    \item $\mathtt{N}_\lambda(\Borel\setminus\{[0,1]\})=\mathtt{N}_\lambda(\Borel\setminus\{(0,1]\})=\mathtt{N}_\lambda(\Borel\setminus\{[0,1)\})=\mathtt{N}_\lambda(\Borel\setminus\{(0,1)\})=\mathtt{N}_\lambda(\Borel)=1$.
    \item $\mathtt{N}_\lambda(\Borel\setminus\{(b,c)\})=c-b$.
    \item $\mathtt{N}_\lambda(\Borel\setminus\{(b,c),[d,e]\})=\min\{c-b, e-d\}$.
\end{enumerate}
\end{example}

\section{Conditional aggregation operators and~generalized survival functions}\label{S-aggr+surv}

From now on, let us consider a collection $\hat{\cE}\subseteq 2^X$ such that $X\in\hat{\cE}$.
Let $\mathbf{G}$ be a~family of real-valued functions on a~set~$X$ containing indicator functions of all sets in~$\hat{\collectionGSF}$. For instance, $\mathbf{G}$ may be the~family of all $\sigma$-measurable functions on~$X$ for a~$\sigma$-algebra~$\Sigma\supset\hat{\collectionGSF}$ on~$X$, like the~setting in~\cite{BoczekHalcinovaHutnikKaluszka2020}.

\begin{definicia}
A~\textit{family of conditional aggregation operators} (FCA for short) is a~family 
$$\bA=\{\aA{\cdot}: E^c\in \hat{\collectionGSF}\}\footnote{
{$\cA$ is a~family of operators parametrized by a~set from~${\hat{\cE}}$.}},$$
such that each $\aA{\cdot}$ is a~map $\aA{\cdot}\colon \mathbf{G}\to[0,+\infty]$
satisfying the following conditions: 
\begin{itemize}
\item[\rm(C1)] $\aA[E]{f}\le \aA[E]{g}$ for any  $f,g\in \mathbf{G}$ such that $f(x)\le g(x)$ for all $x\in E$, $E\neq\emptyset$;  
\item[\rm(C2)]  $\aA{\bin_{E^c}}=0$, $E\neq\emptyset$. 
\end{itemize}
If $E\neq\emptyset$, then $\aA{\cdot}$ is called the~\textit{conditional aggregation operator w.r.t. $E$}
\end{definicia}

\begin{remark}
Moreover, we consider that each element of FCA satisfies $\aA[\emptyset]{\cdot}=0$, unless is stated otherwise (such situation appears only in Section~\ref{Duality}).
\end{remark}

The value $\aA[E]{f}$ can be interpreted as “an aggregated value of $f$ on $E$”. Since the value of $f$ depends on the set $E$, the set $E$ has been called conditional in~\cite{BoczekHalcinovaHutnikKaluszka2020}. Note, the~condition (C1) requires monotonicity of~the~conditional aggregation only on the conditional set (not on the whole space $X$). The~condition (C2) reflects the fact that a zero-valued data on~the~conditional set does not change the~value of~the~conditional aggregation operator. 
As a~consequence of~definition we get $\aA[E]{f} = \aA[E]{f\bin_E}$ for~any $f, {f\bin_{E}} \in \mathbf{G}$ with~a~fixed nonempty $E$.

\bigskip

\begin{example}
Let $\{X\}\subseteq\hat{\cE}\subseteq 2^X$, $\mathbf{G}$ be a family of arbitrary real-valued functions on a set $X$ and $f\in\mathbf{G}$. Then typical examples of the FCA are:\begin{enumerate}[a)]
\item $\cA^{\mathrm{sup}}=\{\aAi[E]{\cdot}{\mathrm{sup}}:E^c\in\hat{\cE}\}$ with $\aAi[E]{f}{\mathrm{sup}}=\sup_{x\in E}f(x)$ for $E\neq\emptyset$;
\item $\cA^{\mathrm{inf}}=\{\aAi[E]{\cdot}{\mathrm{inf}}:E^c\in\hat{\cE}\}$ with $\aAi[E]{f}{\mathrm{inf}}=\inf_{x\in E}f(x)$ for $E\neq\emptyset$;
\item $\cA^{\mathrm{sum}}=\{\aAi[E]{\cdot}{\mathrm{sum}}:E^c\in\hat{\cE}\}$ with $\aAi[E]{f}{\mathrm{sum}}=\sum_{x\in E}f(x)$ for $E\neq\emptyset$ and $X$ being countable.
\end{enumerate}
\end{example}

\begin{example}
Let $\{X\}\subseteq\hat{\cE}\subseteq 2^X$, $\mathbf{G}$ be a family of arbitrary $\sigma$-measurable functions on $X$ for a $\sigma$-algebra $\Sigma\supseteq \hat{\cE}$ on $X$ and $f\in\mathbf{G}$. Then an example of the FCA is
$$\cA^{\mathrm{Ch}_m}=\{\aAi[E]{\cdot}{\mathrm{Ch}_m}:E^c\in\hat{\cE}\}\,\,\,\text{with}\,\,\,\aAi[E]{f}{\mathrm{Ch}_m}=\mathrm{Ch}(f\bin_E,m)\,\,\,\text{for}\,\,\, E\neq\emptyset$$ which is a Choquet integral of the function $f\bin_E$ with respect to the monotone measure $m\colon\Sigma\to[0,+\infty]$. 

Similarly, putting $\aAi[E]{f}{\mathrm{Su}_m}=\mathrm{Su}(f\bin_E,m)$ and $\aAi[E]{f}{\mathrm{Sh}_m}=\mathrm{Sh}(f\bin_E,m)$ (which are Sugeno and Shilkret integral of $f\bin_E$ with respect to $m$), we have $\cA^{\mathrm{Su}_m}$ and $\cA^{\mathrm{Sh}_m}$, respectively. 
\end{example}

It is worth noting that the FCA does not have to contain elements of only one type.

\begin{example}
Let $X=[3]$, $\hat{\cE}=\{\emptyset, \{1\},\{2\},\{3\},X\}$, $\mathbf{G}$ be a family of arbitrary real-valued functions on $X$. We can consider the FCA $$\cA=\{\aAi[E]{\cdot}{\mathrm{sup}}:E^c\in\{\emptyset,\{1\},\{2\}\}\cup\{\aAi[E]{\cdot}{\mathrm{inf}}:E^c\in\{\{3\}\}\}\cup \{\aA[\emptyset]{\cdot}\}.$$
\end{example}

\bigskip
For our following purpose, pick $f\in\mathbf{G}$, and let us set \begin{align}\label{mnozinovafcia}T_{\cA,f}(E):=\aA[E^c]{f}\end{align}for any 
 $E\in\hat{\cE}$.
In fact, $T_{\cA,f}: \hat{\cE}\subseteq 2^X\to[0,+\infty]$. Naturally, $T_{\cA^\mathrm{sup},f}(E)=\sup\limits_{x\in E^c} f(x)$, $T_{\cA^\mathrm{inf},f}(E)=\inf\limits_{x\in E^c} f(x)$, etc.
The FCA is given as $$\bA=\{T_{\cA,\cdot}(E): E\in \hat{\collectionGSF}\}.$$

In~\cite{BoczekHalcinovaHutnikKaluszka2020} the conditional aggregation operators have became the~essence to introduce a~novel concept of generalized survival functions. In the following, we use approach via the set function $T_{\cA,f}$ and then the \textit{generalized survival function} w.r.t. $\cA$ may be defined as follows 
{
\small
\[ 
\mu_\bA(f,\alpha):=
\inf\big\{\mu (E):T_{\cA,f}(E)\le \alpha,\, E\in \hat{\cE}\big\}\footnote{{In the original paper~\cite{BoczekHalcinovaHutnikKaluszka2020} the generalized survival function w.r.t. $\cA$ is defined as $\gsf{}{\x}{\alpha}=\inf\{\mu (E^c):\aA[E]{f}\le \alpha,\, E\in {\cE}\}$, where $\cE=\{E:E^c\in\hat{\cE}\}$ and $\cA=\{\aA[E]{\cdot}:E\in\cE\}.$}},
\]
} 
where $\mu\in\mathbf{M}$ and $f\in\mathbf{G}$, and $\alpha\ge 0$. 
The function $\mu_\bA(f,\alpha)$ is well defined 
as the set $\{E\in \hat{\cE}:T_{\cA,f}(E)\le \alpha\big\}$ is nonempty for all $\alpha\geq 0$, since we suppose $T_{\cA,f}(X)= 0$ and $X\in\hat{\cE}$. Note that standard survival function is a~special case of generalized survival function since $\mu_{\cA^\mathrm{sup}}(f,\alpha)=\mu(\{f>\alpha\})$ whenever $\hat{\cE}$ contains each~$\{f>\alpha\}$.

\textit{A short look at a motivation for studying measure on~hyperset:\,\,} In fact, in the process of computation of generalized survival function we ``measure'' the family of sets. Let us consider the following example adopted to finite space, i.e. $X=\{1,\dots, n\}$ . Then, the input function is a vector ($n$ dimensional), in fact. So, in the following we will use the~notation $\mathbf{x}:=(x_1, \dots, x_n)$. 

\begin{example} Let $X=\{1,2,3\}$, and $\cA^{\mathrm{sum}}=\{\aAi[E^c]{\cdot}{\mathrm{sum}}:E\in\{\emptyset, \{1\}, \{2\}, \{1,2,3\}\}\}$ be FCA.
If $\x=(2,3,4)$, then 
\[
 T_{\cA^{\mathrm{sum}},\x}(E) = 
  \begin{cases} 
   9 & \text{if } E=\emptyset \\
   7       & \text{if } E=\{1\}\\
   6       & \text{if } E={\{2\}}\\
   0       & \text{if } E=\{1,2,3\}.
  \end{cases}
\]
When calculating the generalized survival function for fixed $\alpha$, we are working with the sets of sets, e.g. for $\alpha=7.5$ we are dealing with the set $\hE=\{\{1\}, \{2\}, \{1,2,3\}\}$ in the~following way
\begin{align*}
\mu_{\cA^{\mathrm{sum}}}(\x,7.5)&=\min\big\{\mu (E):T_{\cA^\mathrm{sum} ,\x}(E)\leq 7.5,\, E\in \hat{\cE}\big\}=\min\big\{\mu(\{1\}), \mu(\{2\}),\mu(\{1,2,3\})\}.
\end{align*}
\end{example}

\begin{center}
\begin{figure}[h]
    \centering
    \includegraphics[scale=0.8]{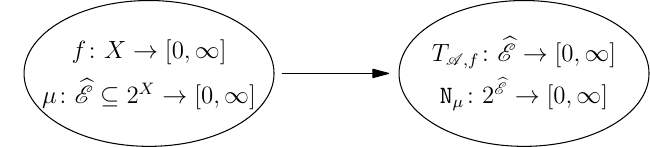}
    \caption{The scheme of transformation of $f$ and $\mu$ into hyperspace}
    \label{fig:scheme}
\end{figure}
\end{center}

\textit{Representation theorem:}\,\,
At the end of this part we provide a result which is the core of~this article. For clarification see the scheme of transformation of $f\in\mathbf{G}$ and $\mu\in\mathbf{M}$ into hyperspace in~Figure~\ref{fig:scheme}.

\begin{veta}[representation theorem]\label{Vtransform}
Let $\cA$ be a FCA with $\hat{\cE}\supseteq\{X\}$, $\mu\in\mathbf{M}$ and $f\in\mathbf{G}$. Then for the generalized survival function w.r.t. $\cA$ it holds that
\[
\gsf{}{f}{\alpha}=\mathtt{N}_\mu(\{E\in\hat{\cE}:T_{\cA,f}(E)>\alpha\big\}).
\]
\end{veta}
\dokaz Immediately we have
{{\begin{align*}
\mathtt{N}_\mu(\{E\in\hat{\cE}:T_{\cA,f}(E)>\alpha\big\})&=\inf\left\{\mu(E)\colon E\in \hat{\cE}\setminus \{T_{\cA,f}>\alpha\big\}\right\}\\
&=\inf\left\{\mu(E):T_{\cA,f}(E)\le \alpha,\, E\in \hat{\cE}\right\}
=\gsf{}{f}{\alpha}, 
\end{align*}}}
which proofs the~equality in the~statement of the~theorem.
\qed

{\begin{remark}In particular, considering size-based super level measures as of~\cite{DoThiele2015,HalcinovaHutnikKiselakSupina2019}, Theorem~\ref{Vtransform} states that 
\[
\mu(s(f)>\alpha)=\mathtt{N}_\mu(\{T_{\cA^\mathrm{size},f}>\alpha\big\}),
\]
where $T_{\cA^\mathrm{size},f}=\sup_{D\in\mathcal{D}} s(f1_{E^c})(D)$
for any 
    $f\in\mathbf{G}$ where $\mathbf{G}$ is a set of all Borel-measurable functions $f\colon X\to[0,+\infty]$ with $\mathcal{D}\subseteq\Borel=\hat{\cE}$. It is called the~outer essential supremum of~$f$ over~$E^c$, where $s$~is a~size. For more details we recommend~\cite{HalcinovaHutnikKiselakSupina2019}.
\end{remark}.}

\section{Generalized Choquet integral}\label{Choquet_integral}

The~generalized survival function defined in the~previous section naturally gives rise to the~generalized Choquet integral defined as follows. 

\begin{definicia}\label{def: Chintegral}\rm(cf.~\cite[Definition 5.4.]{BoczekHalcinovaHutnikKaluszka2020})
 Let $\cA$ be a FCA with $\hat{\cE}\supseteq\{X\}$, $\mu\in\mathbf{M}$ and $f\in {\mathbf{G}}$. The Choquet integral with respect to $\cA$ and $\mu$ ($\cA$-Choquet integral, for short) of $f$ is defined as \begin{align}\label{Chint}\mathrm{C}_{\cA}(f,\mu)=\int_{0}^{\infty}\gsf{}{f}{\alpha}\,\mathrm{d}\alpha.\end{align}\end{definicia}
In the following assertion we shall show that $\mathrm{C}_{\cA}(f,\mu)$ is just the~classical Choquet integral $\mathrm{Ch}(T_{\cA,f},\mathtt{N}_{\mu})$ on the~universe~$\hat{\cE}$, to recall its definition, see~Introduction.
From now on, we require $\emptyset\in\hat{\cE}$ (apart from $X\in\hat{\cE}$), otherwise the value of~the~introduced integral will be infinite.

\begin{dosledok}\label{general=classic}
Let $\cA$ be FCA with $\hat{\cE}\supseteq\{\emptyset,X\}$, $\mu\in\mathbf{M}$ and  $f\in{\mathbf{G}}$. Then $$\mathrm{C}_{\cA}(f,\mu)=\mathrm{Ch}(T_{\cA,f},\mathtt{N}_{\mu}),$$ where the~left-hand side integral is computed on~$X$, and the~right-hand side on $\hat{\cE}\subseteq 2^X$.
\end{dosledok}

\proof We immediately have
{\begin{align*}
\mathrm{C}_{\cA}(f,\mu)=\int_{0}^{\infty}\gsf{}{f}{\alpha}\,\mathrm{d}\alpha=\int_{0}^{\infty}\mathtt{N}_\mu(\{ T_{\cA,f}(E)>\alpha\big\})\,\mathrm{d}\alpha
=\int_{\hat{\cE}}T_{\cA,f}\,\mathrm{d}\mathtt{N}_{\mu}=\mathrm{Ch}(T_{\cA,f},\mathtt{N}_{\mu}).
\end{align*}}\qed

\begin{remark}\label{meratelnost}\rm 
Let us note that in Corollary~\ref{general=classic} we do not have to set the special condition for the measurability of~$T_{\cA,f}$. It follows from~the~fact that the~monotone measure $\mathtt{N}_{\mu}$ is defined on~the~whole power set $2^{\hat{\cE}}$.
On the other hand, if we restrict ourselves to~a~measurable space $(\hat{\cE},\hat{\cS})$, where $\hat{\cS}\subseteq 2^{\hat{\cE}}$ with~$\emptyset\in\hat{\cS}$, then the~previous result holds for~$\hat{\cS}$-measurable set functions $T_{\cA,f}$.
\end{remark}

\begin{remark}\label{size-integral}\rm 
We again stress that a~particular case of~Corollary~\ref{general=classic} is size-based Choquet integral as~of~\cite{DoThiele2015,HalcinovaHutnikKiselakSupina2019}. Thus we have 
\[
\int_{0}^{\infty}\mu(s(f)>\alpha)\,\mathrm{d}\alpha=\mathrm{Ch}(T_{\cA^\mathrm{size},f},\mathtt{N}_{\mu}).
\]
\end{remark}

Now let us consider that $\hat{\cE}\supseteq\{\emptyset, X\}$ is finite. Let us suppose that $|\hat{\cE}|=p$ and {$T_{\cA,\cdot}$ be a~set function with finite values}, and let us assign a~number $i\in\{0,1,\dots,p-1\}$ to~each set $E\in\hat{\cE}$ in~order to hold 
{
\small
\begin{align}\label{usporiadanie}
0=T_{\cA,f}(E_0)\leq\dots\leq T_{\cA,f}(E_{p-1})<+\infty,
\end{align}
} {for $f\in\mathbf{G}$}. 
Then the grasp of generalized survival functions by means of~$\mathtt{N}_{\mu}$-measures immediately gives us the~expression of~the~$\cA$-Choquet integral as~follows. Let us note that in~\cite[Theorem 5.2]{BasarikBorzovaHalcinovaSupina2023} the authors derived several other formulas of the generalized Choquet integral without using the representation theorem.

\begin{tvrdenie}\label{discreteChoquet} Let $\cA$ be FCA with $\hat{\cE}\supseteq\{\emptyset,X\}$ being finite, $\mu\in\mathbf{M}$ and $f\in{\mathbf{G}}$. 
Then $\mathrm{C}_{\cA}(f,\mu)$ is equal to any of~the~following: 
{
\begin{enumerate}[\rm(i)]
    \item $\sum\limits_{i=1}^{p-1}\left(T_{\cA,f}(E_i)-T_{\cA,f}(E_{i-1})\right) \min\limits_{0\leq j\leq i-1}\mu(E_j)$,
    \item $\sum\limits_{i=1}^{p-1} \left[\min\limits_{0\leq j\leq i-1}\mu(E_j)-\min\limits_{0\leq j\leq {i}}\mu(E_j)\right]T_{\cA,f}(E_i)$,
    \item $\sum\limits_{\hat{A}\subseteq\hat{\cE}} \min\limits_{E_i\in\hat{\cE}\setminus\hat{A}}\mu(E_i)\, g_{\hat{A}}(T_{\cA,f})$,
\end{enumerate}
}
where sets $E_i$, $i\in\{0,1,\dots,p-1\}$ satisfy~(\ref{usporiadanie}), and $g_{\hat{A}}(T_{\cA,f})=\max\{0, \min\limits_{E_i\in\hat{A}}T_{\cA,f}(E_i)-\max\limits_{E_i\in\hat{\cE}\setminus\hat{A}}T_{\cA,f}(E_i)\}$\footnote{The maximum and minimum over an empty set are taken as $0$.}. 
\end{tvrdenie}

\proof From Theorem~\ref{Vtransform} and from formulas for standard discrete Choquet integral, see~\cite{BeliakovJamesWu2020,JinKalinaMesiarBorkotokey2018}, denoting $\hat{E_i}=\{E_i,\dots,E_p\}$, it holds: 
{
\begin{align*}
\mathrm{C}_{\cA}(f,\mu)&=\sum_{i=1}^{p-1}\left(T_{\cA,f}(E_i)-T_{\cA,f}(E_{i-1})\right) \mathtt{N}_{\mu}(\hat{E_i})
=\sum_{i=1}^{p-1} \left(\mathtt{N}_{\mu}(\hat{E_i})-\mathtt{N}_{\mu}(\hat{E}_{i+1})\right)T_{\cA,f}(E_i),\\
&=\sum_{\hat{A}\subseteq\hat{\cE}} \mathtt{N}_{\mu}(\hat{A})\, g_{\hat{A}}(T_{\cA,f}).
\end{align*}
}
Using definition of $\mathtt{N}_{\mu}$ given in~(\ref{def:Nmu}) we obtain the result. \qed

\begin{remark}\label{fusion}
The representation of $\cA$-Choquet integral as stated in the previous proposition offers further extensions of this concept. Similar to generalizations of the Choquet integral, see~\cite{DimuroFernandezBedregal2020, LuccaDimuro2019}, the $\cA$-Choquet integral can be generalized considering a pair of fusion functions $H_1, H_2$ such that $H_1$ dominates $H_2$, $H_1$ is $(1,0)$-increasing, $\mu$ be a normalized monotone measure, and $\mathtt{N}_{\mu}$ be symmetric as follows 
{
\begin{align*}
\mathfrak{C}_{\cA}^{H_1, H_2}(f, \mu)&=
\min\{1, \sum_{i=1}^{p-1}H_1(T(E_i),\mathtt{N}_{\mu}(\hat{E}_{(i)}))-H_2(T(E_{i-1}),\mathtt{N}_{\mu}(\hat{E}_{(i)}))\}\\
&=\min\{1, \sum_{i=1}^{p-1}H_1(T(E_i),\min_{j\in\{0,1,\dots, i-1\}}\mu(E_j)-H_2(T(E_{i-1}),\min_{j\in\{0,1,\dots, i-1\}}\mu(E_j))\}
.
\end{align*}
}with $T=T_{\cA,f}$, $\hat{E}_{(i)}=\{E_i,\dots,E_{p-1}\}$ and further settings from Proposition~\ref{discreteChoquet}. Immediately it can be seen that $\mathfrak{C}_{\cA}^{H_1, H_2}(f, \mu)=\mathfrak{C}^{H_1, H_2}(f, \mu)$ with $\cA^{\mathrm{sup}}=\{\sA^{\mathrm{sup}}(\cdot|E): E\in2^X\}$, $X$ being a~finite set. Following this pattern also $g\mathfrak{C}_{\cA}^{H_1, H_2}$ integral can be defined.
\end{remark}

In the previous proposition the representation of $\cA$-Choquet integral via M{\"o}bius transform is omitted, however, we shall deal with it in Section~\ref{Mobius_sec}. Thus we have covered all well-known representations of discrete Choquet integral, see~\cite{BeliakovJamesWu2020}. In the following we aim to underline the usefulness of the above representation from Corollary~\ref{general=classic} of the $\cA$-Choquet integral. We summarize its properties that are simply derived from the properties of the famous Choquet integral. 
From now on, we consider set functions $T_{\cA,\cdot}$ with finite values.
The family of conditional aggregation operators with finite values we shall call $[0,\infty)$-valued FCA. 
Recall that two functions $G, H: \Omega\to[0,\infty)$ are called \textit{comonotone} iff $$(G(s)-G(t))(H(s)-H(t))\geq 0$$ for any $s,t\in\Omega$. A monotone measure $m:\emptyset\in\cS\to[0,\infty]$ will be called \textit{supermodular} iff $m(E\cup F)+m(E\cap F)\geq m(E)+m(F)$ for any $E, F\in\mathcal{S}$. Analogously, a~monotone measure~$m$ will be called \textit{submodular} iff $m(E\cup F)+m(E\cap F)\leq m(E)+m(F)$ for any $E, F\in\mathcal{S}$. A function $H$ is zero {$m$}-almost everywhere iff {$m(\{H>0\})=0$}.

\begin{tvrdenie}\label{vlastnostiC}
Let $\cA,\tilde{\cA}$ 
be $[0,\infty)$-valued FCA with $\hat{\cE}\supseteq\{\emptyset,X\}$. Let $\mu\in\mathbf{M}$, and $f,g\in \mathbf{G}$. 

\begin{enumerate}[\rm(i)]
    \item If $T_{\cA,f}=0$ $\mathtt{N}_{\mu}-$ almost everywhere, then $\Chint=0$. Conversely, if $\mathtt{N}_{\mu}$ is continuous and $\Chint=0$, then $T_{\cA,f}=0$ $\mathtt{N}_{\mu}-$ almost everywhere. 
    \item $\mathtt{N}_{\mu}$ is supermodular if and only if  
    $\mathrm{Ch}(T_{\cA,f}+T_{\tilde{\cA},g},\, \mathtt{N}_{\mu})\geq\mathrm{Ch}(T_{\cA,f},\mathtt{N}_{\mu})+\mathrm{Ch}(T_{\tilde{\cA},g},\mathtt{N}_{\mu})$ for all $T_{\cA,f}, T_{\tilde{\cA},g}$. $\mathtt{N}_{\mu}$ is submodular if and only if  
    $\mathrm{Ch}(T_{\cA,f}+T_{\tilde{\cA},g},\, \mathtt{N}_{\mu})\leq\mathrm{Ch}(T_{\cA,f},\mathtt{N}_{\mu})+\mathrm{Ch}(T_{\tilde{\cA},g},\mathtt{N}_{\mu})$ for all $T_{\cA,f}, T_{\tilde{\cA},g}$.
   \item If $T_{\cA,f}$ and $T_{\tilde{\cA},g}$ are comonotone, then 
    $\mathrm{Ch}(T_{\cA,f}+T_{\tilde{\cA},g},\, \mathtt{N}_{\mu})=\mathrm{Ch}(T_{\cA,f},\mathtt{N}_{\mu})+\mathrm{Ch}(T_{\tilde{\cA},g},\mathtt{N}_{\mu})$.
    \item $\mathrm{C}_{\cA}(f,\mu)=\int_{0}^{\infty}\mathtt{N}_\mu(\{ T_{\cA,f}>\alpha\big\})\,\mathrm{d}\alpha=\int_{0}^{\infty}\mathtt{N}_\mu(\{ T_{\cA,f}\geq\alpha\big\})\,\mathrm{d}\alpha.$
\end{enumerate}
\end{tvrdenie}
\proof All parts follow from~the~standard properties of~the~Choquet integral. Moreover, part ${\rm{(i)}}$ follows from~the~continuity from above of $\mathtt{N}_{\mu}$, see Proposition~\ref{necessvlastnosti}.\qed
\bigskip

As we can see, the properties of the generalized Choquet integral are strongly determined by the properties of monotone measure $\mathtt{N}_{\mu}$ and function $T_{\cA,f}$. Thus, results of Section~\ref{measure_to_measure} play a crucial role. From Proposition \ref{additivitymmi1} we immediately have the following result. It is known that the Choquet integral is linear for each pair of measurable functions $f,g$ if and only if $\mu$ is additive. 

\begin{tvrdenie}\label{linearity}
Let $\cA,\tilde{\cA}$ 
be $[0,\infty)$-valued FCA with $\hat{\cE}\supseteq\{\emptyset,X\}$. Let $\mu\in\mathbf{M}$ with $\mu(X)=c>0$ and  
{$f,g\in\mathbf{G}$.}
Then the~following assertions are equivalent:
\begin{itemize}
\item[\rm (i)]$\mu$ is the strongest monotone measure, i.e., $\mu=\mu^*$. \item[\rm (ii)]$\mathtt{N}_\mu$ is additive.
\item[\rm (iii)] $\mathrm{Ch}(T_{\cA,f}+T_{\tilde{\cA},g},\, \mathtt{N}_{\mu})=\mathrm{Ch}(T_{\cA,f},\mathtt{N}_{\mu})+\mathrm{Ch}(T_{\tilde{\cA},g},\mathtt{N}_{\mu})$ for all $T_{\cA,f}, T_{\tilde{\cA},g}$.
\end{itemize}
\end{tvrdenie}
\proof The equivalence between $\rm (ii)$ and $\rm (iii)$ follows from Proposition~\ref{vlastnostiC} $\rm (ii)$ because supermodularity together with submodularity of monotone measure is equivalent to additivity. The equivalency $\rm (ii)\Leftrightarrow\rm (iii)$ follows from Proposition~\ref{additivitymmi1}.\qed
\bigskip

\section{M{\"o}bius transform}\label{Mobius_sec}

The M{\"o}bius transfrom (or M{\"o}bius inverse) of a monotone measure is widely used in combinatorics since the work~\cite{Rota1964OnTF}. In the framework of discrete monotone measures the M{\"o}bius transform (or M{\"o}bius inverse) is well-known. 
It becomes an~important tool in~game theory (Harsanyi dividends), Dempster-–Shafer theory of~evidence (basic probability assignment), utility theory, etc. It is important from~theoretical point of~view, as well. The M{\"o}bius transform allows introduction of~an~integral w.r.t. discrete monotone measures, which turns to be the~Choquet integral. In~order to use similar benefits of~the~M{\"o}bius transform on general measurable space (not finite) its modification is needed.

The crucial idea on which the results in this part are based, lies in the Corollary~\ref{general=classic}. As the monotone measure $\mathtt{N}_{\mu}$ act on measurable space $(\hat{\cE},2^{\hat{\cE}})$ which need not be finite, we shall adopt the generalized M{\"o}bius transform introduced by Shafer~\cite{Shafer1979} and studied in~\cite{Mesiar2002}. For more details about the standard M{\"o}bius transform of monotone measures we recommend~\cite{Denneberg1994}.
Following the pattern in~\cite{Mesiar2002}, we show that the generalized Choquet integral may be represented as~the~Lebesgue integral by means of the generalized M{\"o}bius transform of measure $\mathtt{N}_\mu$. 
At first, we shall recall the definition of Shafer's allocation between two measurable spaces related to a monotone measure $m$ and M{\"o}bius transform $M$ of $m$, respectively.  

\begin{definicia}
Let  $(\Omega, \mathcal{S})$ and $(Z,\mathcal{Z})$ be  measurable spaces. Let $m: \cS\to[0,\infty)$ be a monotone measure and $M:\mathcal{Z}\to\mathbb{R}$ be a signed measure\footnote{$M:\mathcal{Z}\to\mathbb{R}$ is a signed measure if $\mathcal{Z}$ is a $\sigma$-algebra, $M(\emptyset)=0$, and $M$ is $\sigma$-additive, 
for more details see~\cite{halmos1976measure}.}. If there is a monotone continuous mapping  $h:\mathcal{S}\to \mathcal{Z}$ with $h(\emptyset)=\emptyset$, $h(\Omega)=Z$ such that for all $E\in\mathcal{S}$ it holds that 
\begin{align}\label{allocationvztah}m(E)=M(h(E)),\end{align}
then $h$ is called an \textit{allocation} of $m$ and $M$ is called a~\textit{generalized  M{\"o}bius transform}.
\end{definicia}

Let us consider the measurable space $(\hat{\cE}, \hat{\mathcal{S}})$ with $\hat{\mathcal{S}}\subseteq2^{\hat{\cE}}$, and the~monotone measure~$\mathtt{N}_{\mu}$. Furthermore, let us consider the~measurable space $(Z,\mathcal{Z})$ with the~generalized M{\"o}bius transform $\mathtt{M}$ of $\mathtt{N}_{\mu}$ which is in addition of bounded variation\footnote{I.e., $|M|(X)<\infty$, see~\cite{halmos1976measure}.}. Let us suppose that triples $(\hat{\cE}, \hat{\mathcal{S}}, \mathtt{N}_{\mu})$ and $(Z,\mathcal{Z},\mathtt{M})$ are linked by an~allocation  $h:\hat{\mathcal{S}}\to \mathcal{Z}$. 

For a $\hat{\cS}$-measurable function $T_{\cA,f}:\hat{\cE}\to[0,\infty)$, let us define a mapping $T_{\cA,f}^h:Z\to[0,\infty)$ as 
{
\small
\begin{align}\label{Hh}
T_{\cA,f}^h(z)=\sup\{\alpha: z\in h(\{T_{\cA,f}\geq\alpha\})\}.
\end{align} 
}
Note that  because of the continuity of $h$ for any $\alpha\in(0,\infty)$ it holds that $$\{z\in Z:T_{\cA,f}^h(z)\geq\alpha\}=h(\{T_{\cA,f}\geq\alpha\}),$$ for the proof see Appendix. The previous equality ensures the~$\mathcal{Z}$-measurabitily of~$T_{\cA,f}^h$. Therefore, the~following Lebesgue integral \begin{align}\label{transformint}(\mathcal{L})\,\int_{Z}T_{\cA,f}^h\,\mathrm{d}\mathtt{M}\end{align} is well defined. In~what follows, we show that thanks to~the~existence of~allocation $h$ of~$\mathtt{N}_{\mu}$, the~generalized Choquet integral is just the Lebesgue integral of~$T_{\cA,f}$ with respect to~the~generalized M{\"o}bius transform of $\mathtt{N}_{\mu}$.

\begin{veta}
Let $\cA$ be $[0,\infty)$-valued FCA with $\hat{\cE}\supseteq\{\emptyset,X\}$, $\mu\in\mathbf{M}$ and  $f\in{\mathbf{G}}$. Let $h:\hat{\mathcal{S}}\to\mathcal{Z}$ be an allocation of $\mathtt{N}_{\mu}$ linking triples $(\hat{\cE}, \hat{\mathcal{S}}, \mathtt{N}_{\mu})$ and $(Z,\mathcal{Z},\mathtt{M})$. Then $$\mathrm{C}_{\cA}(f,\mu)=(\mathcal{L})\,\int_{Z}T_{\cA,f}^h\,\mathrm{d}\mathtt{M}.$$
\end{veta}

\proof Using Corollary~\ref{general=classic},~Remark~\ref{meratelnost},~Proposition~\ref{vlastnostiC} (iv) and~(\ref{allocationvztah}) we immediately have 
{\begin{align*}
\mathrm{C}_{\cA}(f,\mu)&=\int_{\hat{\cE}}T_{\cA,f}\,\mathrm{d}\mathtt{N}_{\mu}=\int_0^{\infty}\mathtt{N}_{\mu}(\{T_{\cA,f}(E)>\alpha\})\,\mathrm{d}\alpha=\int_0^{\infty}\mathtt{N}_{\mu}(\{T_{\cA,f}(E)\geq\alpha\})\,\mathrm{d}\alpha\\
&=\int_0^{\infty}\mathtt{M}(h(\{T_{\cA,f}(E)\geq\alpha\}))\,\mathrm{d}\alpha
=\int_0^{\infty}\mathtt{M}(\{T_{\cA,f}^h(z)\geq\alpha\})\,\mathrm{d}\alpha=(\mathcal{L})\,\int_{Z}T_{\cA,f}^h\,\mathrm{d}\mathtt{M}.
\end{align*}}
\qed

Now let us consider $\hat{\cE}$ to be finite. Then $Z$ can be identified with $2^{\hat{\cE}}\setminus\{\emptyset\}$ and an~allocation $h$ of $\mathtt{N}_{\mu}$ is given by $$h(\hat{E})=\{F:\,\emptyset\neq F\subseteq \hat{E}\},$$ see~\cite{Mesiar2002}. Then the formula of $\cA$-Choquet integral by means of M{\"o}bius transform is as follows.

\begin{dosledok}\label{mobius_discrete}
Let $\cA$ be $[0,\infty)$-valued FCA with $\hat{\cE}\supseteq\{\emptyset,X\}$ being finite, $\mu\in\mathbf{M}$ and~$f\in{\mathbf{G}}$.  Let $h:\hat{\mathcal{S}}\to\mathcal{Z}$ be an allocation of $\mathtt{N}_{\mu}$ linking triples $(\hat{\cE}, \hat{\mathcal{S}}, \mathtt{N}_{\mu})$ and $(Z,\mathcal{Z},\mathtt{M})$. Then
\begin{align}\label{Mobiusfin}
\mathrm{C}_{\cA}(f,\mu)=\sum_{\emptyset\neq\hat{E}\subseteq\hat{\cE}}\mathtt{M}(\{\hat{E}\})\min_{E\in\hat{E}}\aA[E^c]{f}.
\end{align}
\end{dosledok}

\proof
From~(\ref{Hh}) for all $T_{\cA,f}:\hat{\cE}\to[0,\infty)$ and~$\hat{E}\subseteq \hat{\cE}$ we obtain {\begin{align*}
T_{\cA,f}^h(\{\hat{E}\})&=\max\{\alpha:\{\hat{E}\}\in h(\{T_{\cA,f}\geq\alpha\})\} =\max\{\alpha:\emptyset\neq\{\hat{E}\}\subseteq \{T_{\cA,f}\geq\alpha\}\}\\
&=\max\left\{\alpha:\left(\forall E\in\hat{E}\right)T_{\cA,f}(E)\geq\alpha\right\}=\min\{T_{\cA,f}(E):\,E\in \hat{E}\}.
\end{align*}}
 Consequently, we get 
 {\begin{align*}
 \mathrm{C}_{\cA}(f,\mu)&=(\mathcal{L})\,\int_{Z}T_{\cA,f}^h\,\mathrm{d}M=\sum_{\emptyset\neq\hat{E}\subseteq\hat{\cE}}T_{\cA,f}^h(\{\hat{E}\})\mathtt{M}(\{\hat{E}\})\\
 &=\sum_{\emptyset\neq\hat{E}\subseteq\hat{\cE}}\mathtt{M}(\{\hat{E}\})\min_{E\in\hat{E}}T_{\cA,f}(E)=
 \sum_{\emptyset\neq\hat{E}\subseteq\hat{\cE}}\mathtt{M}(\{\hat{E}\})\min_{E\in\hat{E}}\aA[E^c]{f}.
 \end{align*}}
 \qed
 
 Note that the classical M{\"o}bius transform of discrete monotone set function $\mathtt{N}_{\mu}$ is given by \begin{align}\label{mobius}\mathtt{M}(\{\hat{E}\})=\sum_{G\subseteq \hat{E}}(-1)^{|\hat{E}\setminus G|}\, \mathrm{\mathtt{N}_{\mu}}(G).\end{align} In fact, the expression~(\ref{mobius}) gives the values of additive set function $\mathtt{M}$ on singletons of $Z$. Recall that the M{\"o}bius transform of $\mathtt{N}_{\mu}$ is invertible by means of the so-called Zeta transform given as \begin{align}\label{vyjadrenie_Nmu}\mathrm{\mathtt{N}_{\mu}}(\hat{E})=\sum_{G\subseteq \hat{E}}\mathtt{M}(\{G\}).\end{align}
 
 Let us remark that formula~(\ref{Mobiusfin}) gives possibilities for further research similarly as in~the~paper~\cite{Fernandez2020}. Here the authors generalized the~standard Choquet integral defined in~terms of~the~M{\"o}bius transform replacing the product by a function. 
 
 The applying of~formula~(\ref{Mobiusfin}) is demostrated in~the~following example.
 
 \begin{example}Let us consider $X=[3]$, and $\cA^{\mathrm{sum}}=\{\aAi[E^c]{\cdot}{\mathrm{sum}}:E\in\{\emptyset, \{1\},[3]\}\}$ be FCA. Let us consider $\x=(1,2,1)$ and the monotone measure $\mu$ with corresponding values in~the~following table.~\begin{table}[h!]
\renewcommand*{\arraystretch}{1.5}
\begin{center}
\begin{tabular}{|c|c|c|}
\hline
$E$ & $\mu(E)$ & $T_{\cA^\mathrm{sum},\x}(E)$\\ \hline
$\emptyset$ & $0$ & $4$\\ \hline
$\{1\}$ & $0.2$ & $3$\\ \hline
$[3]$ & $1$ & $0$\\ \hline
\end{tabular}
\end{center}
\end{table}

Then the generalized survival function takes the~form $$\gsf{\mathrm{sum}}{\x}{\alpha}=\bin_{[0,3)}+0.2\cdot\bin_{[3,4)},$$ and the~corresponding Choquet integral $$\mathrm{C}_{\cA^{\mathrm{sum}}}(\x,\mu)=3\cdot 1+1\cdot 0.2=3.2.$$ Let us compute the generalized Choquet integral by means of~the~M{\"o}bius transform. The corresponding values of minitive measure $\mathtt{N}_{\mu}$, M{\"o}bius transform $\mathtt{M}$ and $\min_{E\in\hat{E}} T_{\cA^\mathrm{sum},\x}(E)=\min_{E\in\hat{E}}\sum_{i\in E^c}x_i$ are given in the below table.\begin{table}[h!]
\renewcommand*{\arraystretch}{1.5}
\begin{center}
\begin{tabular}{|c|c|c|c|}
\hline
$\hat{E}$ & $\mathtt{N}_{\mu}(\hat{E})$ & $\mathtt{M}(\{\hat{E}\})$ & $\min\limits_{E\in\hat{E}} T_{\cA^\mathrm{sum},\x}(E)$\\ \hline
$\{\emptyset\}$ & 0.2 & 0.2 & 4\\ \hline
\{\{1\}\} & 0 & 0 & 3\\ \hline
\{[3]\} & 0 & 0 & 0\\ \hline
$\{\emptyset,\{1\}\}$ & 1 & 0.8 & 3\\ \hline
$\{\emptyset,[3]\}$ & 0.2 & 0 & 0\\ \hline
\{\{1\},[3]\} & 0 & 0 & 0\\ \hline
$\{\emptyset,\{1\},[3]\}$ & 1 & 0 & 0\\ \hline
\end{tabular}
\end{center}
\end{table}
 \end{example}
 Then 
 {
 \small
 $$\sum_{\emptyset\neq\hat{E}\subseteq\hat{\cE}}\mathtt{M}(\{\hat{E}\})\min\limits_{E\in\hat{E}}\sum\limits_{i\in E^c}x_i=4\cdot 0.2+0.8\cdot 3=3.2.$$
 }

\section{Duality}\label{Duality}

We assume thorughout the~whole section that $\hat{\cE}\supseteq\{X\}$ is closed under complements. We shall compute a~dual measure of~$\mathtt{N}_{\mu}$. Let us recall that a~dual measure~$m^d$ of a~measure~$m$ with~$m(\Omega)<\infty$ is defined as 
\begin{align*}
m^d(E)=m(\Omega)-m(E^c).
\end{align*}
We shall show that a maxitive measure defined as \begin{align*}
\mathtt{\Pi}_{\mu^d}(\hE)=\sup\{\mu^d(E^c):E\in\hE\}.
\end{align*}  
is dual to minitive measure $\mathtt{N}_{\mu}$.
Note that another maxitive measure on a~hyperset~$2^{\hat{\cE}}$ was considered in~\cite{BoczekHutnikKaluzskaKleinova2021}.

\begin{lema}\label{dual}
Let $\hat{\cE}\supseteq\{X\}$ 
be closed under complements and $\mu(X)<\infty$. The~maxitive measure~$\mathtt{\Pi}_{\mu^d}$ is a~dual measure to~$\mathtt{N}_{\mu}$.
\end{lema}
\proof
Let us take arbitrary set~$\hE\subseteq\hat{\cE}$. Then we may derive that
{\begin{align*}
\mathtt{N}_\mu(\hE)&=\inf\{\mu(E):E\in\hat{\cE}\setminus\hE\}
=\inf\{\mu(X)-\mu^d(E^c):E\in\hat{\cE}\setminus\hE\}\\ 
&=\mu(X)-\sup\{\mu^d(E^c):E\in\hat{\cE}\setminus\hE\}=\mathtt{N}_\mu(\hat{\cE})-\sup\{\mu^d(E^c):E\in\hat{\cE}\setminus\hE\}.\ 
\end{align*}}
The~latter formula is equal to~$\mathtt{N}_\mu(\hat{\cE})-\mathtt{\Pi}_{\mu^d}(\hat{\cE}\setminus\hE)$, which shows that the~measure~$\mathtt{\Pi}_{\mu^d}$ is dual to~$\mathtt{N}_{\mu}$.
\qed

The~generalized survival function~$\gsf{}{f}{\alpha}$ imitates strict survival function~$\mu(\{f>\alpha\})$. On the~other hand, M.~Boczek, O.~Hutn\'ik, M.~Kaluszka and M.~Kleinov\'a~M.~\cite{BoczekHutnikKaluzskaKleinova2021} introduced a~notion of~$\mu_{\bar{\cA}}^+(f,\alpha)$ which imitates a non-strict survival function~$\mu(\{f\geq\alpha\})$. They have shown that both notions are in certain sense complementary. We reprove their result using measure on hyperset instead of original direct approach.

Let $\overline{y}\in(0,\infty)$. We say that $\cA$ is $\overline{y}$-idempotent FCA if $T_{\cA,\overline{y}}(E)=\overline{y}$ for any $E\in\hat{\cE}\setminus\{\emptyset\}$. We set $\mathbf{G}^{\overline{y}}$ to be a~family of~all nonnegative functions~$f\in\mathbf{G}$ such that $\sup_{x\in X } f(x)\leq \overline{y}$. Note that if $\cA$ is $\overline{y}$-idempotent FCA, then $T_{\cA,\cdot}$ is $[0,\overline{y}]$-valued set function on~{$\mathbf{G}^{\overline{y}}$}. The~second step is considering a~conditional aggregation operator w.r.t. $E\neq\emptyset$ given as $$\overline{\sA}(f|E):=\sA(\overline{y}|E)-\sA(\overline{y}-f|E),$$ and $\overline{\sA}(f|\emptyset)=\infty$. The set of such conditional aggregation operators we shall denote by $$\bar{\cA}=\{\overline{\sA}(\cdot|E):E\in\hat{\cE}\},$$ where $\hat{\cE}\supseteq\{X\}$ is an~arbitrary collection of sets which is closed under complements. 
The corresponding transformation we shall denote by
\begin{align*}
\overline{T}_{\bar{\cA},f}(E):=\overline{\sA}(f|E)
\end{align*} 
for any $E\in\hat{\cE}$. It can be easily seen that the set function $\overline{T}_{\bar{\cA},f}(E)$ satisfies the properties (C1), (C2) of conditional aggregation operators. In keeping with our notion, the~generalized level measure ${\nu_{\bar{\cA}}^+(f,\alpha)}$ defined in~\cite{BoczekHalcinovaHutnikKaluszka2020} may be given by the~formula 
{
\small
\begin{align}\label{glm}
{\nu_{\bar{\cA}}^+(f,\alpha)}=\sup\{\nu(E):\overline{T}_{\bar{\cA},f}(E)\geq \alpha, E\in\hat{\cE}\},
\end{align} 
}
where $\alpha\geq 0$.

\begin{tvrdenie}\label{dual-levels}
Let $\cA$ be $\overline{y}$-idempotent FCA with $\hat{\cE}\supseteq\{X\}$ being closed under complements, $\mu(X)<\infty$, and $f\in\mathbf{G}^{\overline{y}}$. Then
$$\gsf{}{f}{\alpha}+\left(\mu^d\right)_{\bar{\cA}}^+(\overline{y}-f,\overline{y}-\alpha)=\mu(X)$$
for any $\alpha\in[0,\overline{y}]$.

\end{tvrdenie}
\proof
Using the~dual measure $\mathtt{\Pi}_{\mu^d}$, by~Theorem~\ref{Vtransform}, and Lemma~\ref{dual} we obtain
{\begin{align}\label{vztah}
\gsf{}{f}{\alpha}&=\mathtt{N}_\mu \left(\{T_{\cA,f}(E)>\alpha\}\right)\nonumber
=\mathtt{N}_\mu(\hat{\cE})-\mathtt{\Pi}_{\mu^d}\left(\hat{\cE}\setminus\{T_{\cA,f}(E)>\alpha\}\right)\nonumber\\
&=\mathtt{N}_\mu(\hat{\cE})-\mathtt{\Pi}_{\mu^d}\left(\{T_{\cA,f}(E)\leq\alpha\}\right).\
\end{align}}
One can check that $\overline{T}_{\bar{\cA},f}(\emptyset)=\infty$, and 
\begin{align}\label{T}
\overline{T}_{\bar{\cA},f}(E)=T_{\cA,\overline{y}}(E^c)-T_{\cA,\overline{y}-f}(E^c),
\end{align}
for $E\neq\emptyset$. If we take $g(x)=\overline{y}-f(x)$ for all $x\in X$, then from~(\ref{vztah}) we have that $\gsf{}{f}{\alpha}$ is equal to
{
\small
\[
\mu(X)-\mathtt{\Pi}_{\mu^d}\left(\{E\in\hat{\cE}:T_{\cA,\overline{y}-g}(E)\leq\alpha\}\right),
\]
}
and
{
\begin{align}\label{dokazPi}
\mathtt{\Pi}_{\mu^d}\left(\{E\in\hat{\cE}:T_{\cA,\overline{y}-g}(E)\leq\alpha\}\right)&=\mathtt{\Pi}_{\mu^d}\left(\{E\in\hat{\cE}:T_{\cA,\overline{y}}(E)-\overline{T}_{\bar{\cA},g}(E^c)\leq\alpha\}\right)\\
&=\mathtt{\Pi}_{\mu^d}\left(\{E\in\hat{\cE}:\overline{T}_{\bar{\cA},g}(E^c)\geq \overline{y}-\alpha\}\right)\nonumber\\
&=\sup\{\mu^d(E^c):\overline{T}_{\bar{\cA},\overline{y}-f}(E^c)\geq\overline{y}-\alpha,E^c\in\hat{\cE}\}\nonumber\\
&=\left(\mu^d\right)_{\bar{\cA}}^+ (\overline{y}-f,\overline{y}-\alpha).\nonumber\
\end{align}
}
for any $\alpha\in[0,\overline{y}]$.
In the~equality~(\ref{dokazPi}) we used the~relation~(\ref{T}), which holds for $E\neq\emptyset$. However, the~case $E=\emptyset$ is not a~problem, since $\emptyset^c=X\in\{E\in\hat{\cE}:T_{\cA,\overline{y}-g}(E)\leq\alpha\}$ for any $\alpha\in[0,\overline{y}]$ because $T_{\cA,\overline{y}-g}(X)=0$ and $\emptyset^c=X\in\{E\in\hat{\cE}:T_{\cA,\overline{y}}(E)-\overline{T}_{\bar{\cA},g}(E^c)\leq\alpha\}$ for any $\alpha\in[0,\overline{y}]$ because $T_{\cA,\overline{y}}(X)-\overline{T}_{\bar{\cA},g}(\emptyset)=-\infty$.
\qed

There are two major differences between our Proposition~\ref{dual-levels} and corresponding assertion in~\cite{BoczekHutnikKaluzskaKleinova2021}. First of all, the~formula in the~statement of Proposition~\ref{dual-levels} contains the~generalized survival function~$\gsf{}{f}{\alpha}$ for arbitrary function~$f$ and real~$\alpha$. The~concept by~\cite{BoczekHutnikKaluzskaKleinova2021} contains computed parameters so that the~presented equality holds. However, the~corresponding formula in~\cite{BoczekHutnikKaluzskaKleinova2021} contains their concept of~$\mu_{\bar{\cA}}^+(f,\alpha)$  for arbitrary function~$f$ and real~$\alpha$, and the~parameters of the~generalized survival function are computed. The~second difference is the~proof, since we apply measure on~hyperset, but the~corresponding proof in~\cite{BoczekHutnikKaluzskaKleinova2021} is direct and relies on the~original definitions of the~applied notions.   

\section*{Conclusion}
We have shown that the $\cA$-Choquet integral can be represented as a standard Choquet integral of a transformed input function $T_{\cA,f}$ with respect to transformed monotone measure~$\mathtt{N}_{\mu}$. The~crucial results from this point of view are Theorem~\ref{Vtransform} and Corollary~\ref{general=classic}. As in many applications in the past the minitive/necessity monotone measures naturally appeared. In fact, it turned out that the corresponding transformed monotone measure $\mathtt{N}_{\mu}$ is minitive monotone measure.

Thanks to perspective of above mentioned representation we have formulated some (new) properties of this integral, we have represented the generalized Choquet integral via several formulas among them via M{\"o}bius transform, see Proposition~\ref{vlastnostiC}, Proposition~\ref{linearity}, Proposition~\ref{discreteChoquet}, Corollary~\ref{mobius_discrete}. The study of further properties of $\cA$-Choquet integral based on~transformation theorem remains for future research.

\section*{Acknowledgement}
The support of the grants APVV-21-0468 and VEGA 1/0657/22 
is kindly announced.

\bigskip
\bigskip


\begin{thebibliography}{10}

\bibitem{BasarikBorzovaHalcinovaSupina2023}
{\sc Basarik, S., Borzov{\'{a}}, J., Hal{\v{c}}inov{\'{a}}, L., and \v{S}upina,
  J.}
\newblock {Conditional aggregation-based Choquet integral on discrete space},
  submitted.

\bibitem{BeliakovJamesWu2020}
{\sc Beliakov, G., James, S., and Wu, J.-Z.}
\newblock {\em {Discrete fuzzy measures: Computational aspects}}.
\newblock Springer, Cham, 2020.

\bibitem{BoczekHalcinovaHutnikKaluszka2020}
{\sc Boczek, M., Hal{\v{c}}inov{\'a}, L., Hutn{\'i}k, O., and Kaluszka, M.}
\newblock Novel survival functions based on conditional aggregation operators.
\newblock {\em Information Sciences 580\/} (2021), 705--719.

\bibitem{BoczekHutnikKaluszka2022}
{\sc Boczek, M., Hutn\'{i}k, O., and Kaluszka, M.}
\newblock {Choquet-Sugeno-like operator based on relation and conditional
  aggregation operators}.
\newblock {\em Information Sciences 582\/} (2022), 1--21.

\bibitem{BoczekHutnikKaluzskaKleinova2021}
{\sc Boczek, M., Hutn\'{i}k, O., Kaluszka, M., and Kleinov{\'{a}}, M.}
\newblock Generalized level measure based on a family of conditional
  aggregation operators.
\newblock {\em Fuzzy Sets and Systems 457\/} (2023), 180--196.

\bibitem{Denneberg1994}
{\sc Denneberg, D.}
\newblock {\em {Non-Additive Measure and Integral}}.
\newblock {Kluwer Academic Publishers}, Dordrecht, 1994.

\bibitem{DimuroFernandezBedregal2020}
{\sc Dimuro, G.~P., Fern\'{a}ndez, J., Bedregal, B., Mesiar, R., Sanz, J.~A.,
  Lucca, G., and Bustince, H.}
\newblock {The state-of-art of the generalizations of the Choquet integral:
  From aggregation and pre-aggregation to ordered directionally monotone
  functions}.
\newblock {\em Information Fusion {\bf 57}\/} (2020), 27--43.

\bibitem{DoThiele2015}
{\sc Do, Y., and Thiele, C.}
\newblock {$L^p$ theory for outer measures and two themes of Lennart Carleson
  united}.
\newblock {\em Bulletin of the American Mathematical Society 52}, 2 (2015),
  249--296.

\bibitem{Dubois2015}
{\sc Dubois, D., and Prade, H.}
\newblock {\em {Possibility Theory and Its Applications: Where Do We Stand?}}
\newblock Springer Berlin Heidelberg, Berlin, Heidelberg, 2015, pp.~31--60.

\bibitem{Fernandez2020}
{\sc Fernandez, J., Bustince, H., Horanská, L., Mesiar, R., and Stup{\v
  n}anová, A.}
\newblock {A Generalization of the Choquet Integral Defined in Terms of
  the~M{\"o}bius Transform}.
\newblock {\em IEEE Transactions on Fuzzy Systems {\bf 28}}, 10 (2020),
  2313--2319.

\bibitem{Halcinova2017}
{\sc Hal{\v{c}}inov{\'a}, L.}
\newblock Sizes, super level measures and integrals.
\newblock In {\em Aggregation Functions in Theory and in Practice, 9th
  International Summer School on Aggregation Functions, Sk{\"{o}}vde, Sweden,
  19--22 June 2017\/} (2017), V.~Torra, R.~Mesiar, and B.~DeBaets, Eds.,
  vol.~581 of {\em Advances in Intelligent Systems and Computing}, Springer,
  pp.~181--188.

\bibitem{HalcinovaHutnikKiselakSupina2019}
{\sc Hal{\v{c}}inov{\'a}, L., Hutn{\'i}k, O., Kise{\v{l}}{\'a}k, J., and
  {\v{S}}upina, J.}
\newblock Beyond the scope of super level measures.
\newblock {\em Fuzzy Sets and Systems 364\/} (2019), 36--63.

\bibitem{halmos1976measure}
{\sc Halmos, P.}
\newblock {\em {Measure Theory}}.
\newblock Graduate Texts in Mathematics. Springer New York, 1976.

\bibitem{HondaOkazaki2017}
{\sc Honda, A., and Okazaki, Y.}
\newblock Theory of inclusion–exclusion integral.
\newblock {\em Information Sciences \bf 376\/} (2017), 136--147.

\bibitem{JinKalinaMesiarBorkotokey2018}
{\sc Jin, L., Kalina, M., Mesiar, R., and Borkotokey, S.}
\newblock {Discrete Choquet Integrals for Riemann Integrable Inputs With Some
  Applications}.
\newblock {\em IEEE Transactions on Fuzzy Systems {\bf 26}}, 5 (2018),
  3164--3169.

\bibitem{KlementMesiarPap2010}
{\sc Klement, E.~P., Mesiar, R., and Pap, E.}
\newblock {A~universal integral as common frame for Choquet and Sugeno
  integral}.
\newblock {\em IEEE Transactions on Fuzzy Systems 18\/} (2010), 178--187.

\bibitem{KRATSCHMER2003455}
{\sc Krätschmer, V.}
\newblock {When fuzzy measures are upper envelopes of probability measures}.
\newblock {\em Fuzzy Sets and Systems {\bf 138}}, 3 (2003), 455--468.

\bibitem{LuccaDimuro2019}
{\sc Lucca, G., Dimuro, G.~P., Fernández, J., Bustince, H., Bedregal, B., and
  Sanz, J.~A.}
\newblock {Improving the Performance of Fuzzy Rule-Based Classification Systems
  Based on a Nonaveraging Generalization of CC-Integrals Named
  $C_{F_1F_2}$-Integrals}.
\newblock {\em IEEE Transactions on Fuzzy Systems {\bf 27}}, 1 (2019),
  124--134.

\bibitem{Mesiar2002}
{\sc Mesiar, R.}
\newblock {Fuzzy measures and generalized M{\"o}bius transform}.
\newblock {\em International Journal of General Systems {\bf 31}}, 6 (2002),
  587--599.

\bibitem{MesiarSipos1994}
{\sc Mesiar, R., and {\v{S}}ipo{\v{s}}, J.}
\newblock A theory of fuzzy measures: Integration and its additivity.
\newblock {\em International Journal of General Systems {\bf 23}}, 1 (1994),
  49--57.

\bibitem{Pap1995NullAdditiveSF}
{\sc Pap, E.}
\newblock {\em {Null-Additive Set Functions}}.
\newblock Kluwer/Ister Science, Dordrecht-Bratislava, 1995.

\bibitem{Poncet2017}
{\sc Poncet, P.}
\newblock {Representation of maxitive measures: An overview}.
\newblock {\em Mathematica Slovaca {\bf 67}}, 1 (2017), 121--150.

\bibitem{Puhalskii}
{\sc Puhalskii, A.}
\newblock {\em Large Deviations and Idempotent Probability (1st ed.)}.
\newblock Chapman and Hall/CRC, New York, 2001.

\bibitem{Rota1964OnTF}
{\sc Rota, G.}
\newblock {On the foundations of combinatorial theory I. Theory of M{\"o}bius
  Functions}.
\newblock {\em Zeitschrift f{\"u}r Wahrscheinlichkeitstheorie und Verwandte
  Gebiete {\bf 2}\/} (1964), 340--368.

\bibitem{Shafer1979}
{\sc Shafer, G.}
\newblock {Allocations of Probability}.
\newblock {\em The Annals of Probability {\bf 7}}, 5 (1979), 827 -- 839.

\bibitem{Shilkret1971}
{\sc Shilkret, N.}
\newblock Maxitive measure and integration.
\newblock {\em Indagationes Mathematicae {33}\/} (1971), 109--116.

\bibitem{Sipos1979}
{\sc {\v{S}}ipo{\v{s}}, J.}
\newblock Integral with respect to a pre-measure.
\newblock {\em Mathematica Slovaca \bf 29\/} (1979), 141--155.

\bibitem{Sugeno1974}
{\sc Sugeno, M.}
\newblock {\em Theory of fuzzy integrals and its applications}.
\newblock PhD thesis, Tokyo Institute of Technology, Tokyo, Japan, 1974.

\bibitem{WangKlir2009}
{\sc Wang, Z., and Klir, G.}
\newblock {\em {Generalized Measure Theory}}.
\newblock Springer, New York, 2009.

\bibitem{YagerMesiar2015}
{\sc Yager, R.~R., and Mesiar, R.}
\newblock {On the Transformation of Fuzzy Measures to the Power Set and Its
  Role in Determining the Measure of a Measure}.
\newblock {\em IEEE Transactions on Fuzzy Systems {\bf 23}}, 4 (2015),
  842--849.

\bibitem{ZADEH19783}
{\sc Zadeh, L.}
\newblock Fuzzy sets as a basis for a theory of possibility.
\newblock {\em Fuzzy Sets and Systems {\bf 1}}, 1 (1978), 3--28.

\end{thebibliography}

\section*{Appendix}

\dokazp{Proposition~\ref{necesschar}}
(a) Let $m$ be a~minitive measure on~$2^\Omega$, and let us set $\pi'(x)=m(\Omega\setminus\{x\})$ for $x\in\Omega$. Then 
{\small\begin{align*}
m(E)&=m\left(\bigcap\limits_{x\in \Omega\setminus E}(\Omega\setminus\{x\})\right)=\inf\limits_{x\in \Omega\setminus E}m(\Omega\setminus \{x\})\\
&=\inf\limits_{x\in \Omega\setminus E}\pi'(x).
\end{align*}}
\par
We assume now that $\pi'\colon \Omega\to[0,\infty]$, and we set $m(E)=\inf\{\pi'(x): x\in \Omega\setminus E\}$ for any $E\subseteq \Omega$. Then
{\begin{align*}
m\left(\bigcap\limits_{t\in T}E_t\right)=\inf\{\pi'(x): x\in \Omega\setminus \bigcap\limits_{t\in T}E_t\}=\inf\limits_{t\in T}\inf\{\pi'(x): x\in \Omega\setminus E_t\}
=\inf\limits_{t\in T}m(E_t).
\end{align*}}
\par
(b) Let $m$ be a~maxitive measure on~$2^\Omega$, and let us set $\pi(x)=m(\{x\})$ for $x\in\Omega$. Then 
{\begin{align*}
m(E)&=m\left(\bigcup\limits_{x\in E}\{x\}\right)=\sup\limits_{x\in E}m(\{x\})=\sup\limits_{x\in E}\pi(x).
\end{align*}}
\par
We assume now that $\pi\colon \Omega\to[0,\infty]$, and we set $m(E)=\sup\{\pi(x): x\in E\}$ for any $E\subseteq \Omega$. Then
{\begin{align*}
m\left(\bigcup\limits_{t\in T}E_t\right)=\sup\{\pi(x): x\in\bigcup\limits_{t\in T}E_t\}=\sup\limits_{t\in T}\sup\{\pi(x): x\in E_t\}=\sup\limits_{t\in T}m(E_t).
\end{align*}}
\qed

\dokazp{equality $\{z\in Z:T_{\cA,f}^h(z)\geq\alpha\}=h(\{E\in\hat{\cE}:T_{\cA,f}(E)\geq\alpha\})$ appearing below equality~\eqref{Hh}}
One can easily see that if $z\in h(\{T_{\cA,f}\geq\alpha\})$ then 
\[
T_{\cA,f}^h(z)=\sup\{\beta: z\in h(\{T_{\cA,f}\geq\beta\})\}\geq\alpha.
\]
To prove the~reversed inclusion, let us take $z\in Z$ with $T_{\cA,f}^h(z)\geq\alpha$. By the~definition of~$T_{\cA,f}^h(z)$ as supremum, for each $n$ there is $\alpha_n\geq\alpha-\frac{1}{2^n}$ with $z\in h(\{T_{\cA,f}\geq\alpha_n\})$. However, the~continuity of~$h$ yields
\[
\bigcap_nh(\{T_{\cA,f}\geq\alpha_n\})=h(\{T_{\cA,f}\geq\alpha\}).
\]
Thus since $z\in h(\{T_{\cA,f}\geq\alpha_n\})$ for each~$n$ we obtain $z\in h(\{T_{\cA,f}\geq\alpha\})$.
\qed

\bigskip
\bigskip

\vspace{5mm}
\noindent \small{\textsc{Jana Borzov\'a, Lenka Hal\v cinov\'a, Jaroslav \v Supina} \newline Institute of Mathematics, Faculty of
Science, Pavol Jozef \v Saf\'arik University in Ko\v sice,
\newline {\it Current address:} Jesenn\'a 5, SK 040~01 Ko\v sice,
Slovakia,
\newline {\it E-mail addresses:} jana.borzova@upjs.sk \newline
\phantom{{\it E-mail addresses:}} lenka.halcinova@upjs.sk \newline
\phantom{{\it E-mail addresses:}} jaroslav.supina@upjs.sk}

\end{document}